\documentclass[12pt]{amsart}

\usepackage{amsmath,amssymb}

\usepackage[all]{xy}

\numberwithin{equation}{section}
\theoremstyle{plain}
\newtheorem{thm}{Theorem}[section]
\newtheorem{prop}[thm]{Proposition}
\newtheorem{lemm}[thm]{Lemma}
\newtheorem{coro}[thm]{Corollary}

\theoremstyle{definition}
\newtheorem{defi}[thm]{Definition}
\newtheorem{exas}[thm]{Examples} 

\newtheorem{rems}[thm]{Remarks}  
\newtheorem{rem}[thm]{Remark} 

\newcommand{\bd}{\begin{defi}}
\newcommand{\ed}{\end{defi}}
\newcommand{\bt}{\begin{thm}}
\newcommand{\et}{\end{thm}}
\newcommand{\bp}{\begin{prop}}
\newcommand{\ep}{\end{prop}}
\newcommand{\bl}{\begin{lemm}}
\newcommand{\el}{\end{lemm}}
\newcommand{\bc}{\begin{coro}}
\newcommand{\ec}{\end{coro}}
\newcommand{\bpf}{\begin{proof}}
\newcommand{\epf}{\end{proof}}
\newcommand{\br}{\begin{rems}\begin{flushleft}\end{flushleft}\nopagebreak} 
\newcommand{\er}{\end{rems}}
\newcommand{\bex}{\begin{exas}\begin{flushleft}\end{flushleft}\nopagebreak} 
\newcommand{\eex}{\end{exas}}

\newcommand{\be}{\begin{enumerate}}
\newcommand{\ee}{\end{enumerate}}
\newcommand{\bi}{\begin{itemize}}
\newcommand{\ei}{\end{itemize}}

\newcommand{\xyinc}{\ar@{^{(}->}} 

\newcommand{\smallxymatrix}[1]{\xymatrix@1@=1pc {#1} } 

\newcommand{\tinyxymatrix}[1]{
\def\objectstyle{\scriptstyle}
\def\labelstyle{\scriptstyle}
\vcenter{\xymatrix@-1.9pc{#1} }}

\newcommand{\descent}{\tinyxymatrix{*{\circ}\ar@{-}[rd] & &\\
&*{\circ}\ar@{-}[rd] &\\
& & *{\circ} }}
\newcommand{\ascent}{\tinyxymatrix{& & *{\circ}\ar@{-}[ld] \\
&*{\circ}\ar@{-}[ld] &\\
 *{\circ}& & }}
\newcommand{\valley}{\tinyxymatrix{*{\circ}\ar@{-}[rd] & &*{\circ}\\
&*{\circ}\ar@{-}[ru] &}}
\newcommand{\peak}{\tinyxymatrix{&*{\circ}\ar@{-}[rd] &\\
*{\circ}\ar@{-}[ru] & &*{\circ} }}

\newcommand{\sdescent}[1]{\tinyxymatrix{*{\circ}\ar@{-}[rd] & &\\
&*{\circ}\ar@{-}[rd]\save[]+<2pt,5pt>*{#1}\restore &\\
& & *{\circ} }}
\newcommand{\sascent}[1]{\tinyxymatrix{& & *{\circ}\ar@{-}[ld] \\
&*{\circ}\ar@{-}[ld]\save[]+<-2pt,5pt>*{#1}\restore &\\
 *{\circ}& & }}
\newcommand{\svalley}[1]{\tinyxymatrix{*{\circ}\ar@{-}[rd] & &*{\circ}\ar@{-}[ld]\\
&*{\circ}\save[]+<0pt,-6pt>*{#1}\restore &}}
\newcommand{\speak}[1]{\tinyxymatrix{&*{\circ}\ar@{-}[rd]\save[]+<0pt,6pt>*{#1}\restore &\\
*{\circ}\ar@{-}[ru] & &*{\circ} }}

\newcommand{\sumsub}[1]{\sum_{\substack{\rule{0pt}{8pt}#1\rule{0pt}{8pt}}}} 

\newcommand{\ten}{\mbox{\hspace*{-.5pt}\raisebox{1pt}{${\scriptstyle \otimes}$}   
\hspace*{-4pt}}}                                                                  
\newcommand{\Ualpha}                                                           
{\mbox{{$\alpha$}\hspace*{-8.3pt}\raisebox{-5pt}{$\scriptstyle{\smile}$}}}
\newcommand{\Oalpha}
{\mbox{{$\alpha$}\hspace*{-8pt}\raisebox{6.5pt}{$\scriptstyle{\frown}$}}}

\newcommand{\id}{\mathit{id}}          
\newcommand{\abs}[1]{\lvert#1\rvert}   
\newcommand{\sign}{\mathrm{sign}}
\newcommand{\ipart}[1]{\lfloor \frac{#1}{2} \rfloor} 
\newcommand{\ipartn}{\ipart{n}}  
\newcommand{\liebrac}[1]{[#1]}

\newcommand{\Des}{\mathrm{Des}}
\newcommand{\Sh}{\mathrm{Sh}}

\newcommand{\Peak}{\mathrm{Peak}}
\newcommand{\Peakint}{\overset{\circ}{\Peak}}

\newcommand{\Span}{\mathrm{Span}}
\renewcommand{\ker}{\mathrm{Ker}} 

\newcommand{\im}{\mathrm{Im}}

\newcommand{\End}{\textrm{End}}

\newcommand{\Q}{\mathbb{Q}}
\newcommand{\Z}{\mathbb{Z}}

\newcommand{\QS}{\mathbb{Q}\calS{}}
\newcommand{\QB}{\mathbb{Q}\calB{}}
\newcommand{\calS}[1]{\mathcal{S}_{#1}}
\newcommand{\calSn}{\calS{n}}
\newcommand{\calB}[1]{\mathcal{B}_{#1}}
\newcommand{\calBn}{\calB{n}}
\newcommand{\calD}[1]{\mathcal{D}_{#1}}
\newcommand{\calDn}{\calD{n}}
\newcommand{\calQ}{\mathcal{Q}}
\newcommand{\calF}{\mathcal{F}}
\newcommand{\calQn}{\calQ_{n}}
\newcommand{\calFn}{\calF_{n}}

\newcommand{\bfX}{\mathbf{X}}
\newcommand{\bfY}{\mathbf{Y}}
\newcommand{\bfTX}{\mathbf{TX}}
\newcommand{\bfTY}{\mathbf{TY}}
\newcommand{\bfx}{\mathbf{x}}
\newcommand{\bfy}{\mathbf{y}}

\newcommand{\ppp}[1]{\mathfrak{P}_{#1}}
\newcommand{\pppn}{\ppp{n}}
\newcommand{\pppint}[1]{\overset{\circ}{\mathfrak{P}}_{#1}}
\newcommand{\pppnint}{\pppint{n}}
\newcommand{\calFint}[1]{\overset{\circ}{\mathcal{F}_{#1}}}
\newcommand{\calFnint}{\calFint{n}}
\newcommand{\Lambdaint}{\overset{\circ}{\Lambda}}

\newcommand{\ThetaB}{\Theta^{\pm}}
\newcommand{\Stilde}{R} 

\newcommand{\Pint}{\,\overset{\circ}{\!P}}
\newcommand{\pint}{\,\overset{\circ}{\!p}}
\newcommand{\wpint}{\,\overset{\circ}{\!\wp}}
\newcommand{\hwp}{\widehat{\wp}}

\newcommand{\QSym}{\mathcal{Q}\mathit{Sym}}
\newcommand{\SSym}{\mathfrak{S}\mathit{Sym}}

\newcommand{\Sol}[1]{\mathit{Sol}(#1)}
\newcommand{\SolB}{\Sol{B_n}}
\newcommand{\SolA}{\Sol{A_{n-1}}}
\newcommand{\SolD}{\Sol{D_n}}
\newcommand{\Ome}[1]{\mathit{Sol}^{\pm}(#1)}
\newcommand{\OmeB}{\Ome{B_n}}
\newcommand{\sol}[1]{\mathit{s}(#1)}
\newcommand{\sB}{\sol{B_n}}
\newcommand{\sA}{\sol{A_{n-1}}}
\newcommand{\sD}{\sol{D_n}}
\newcommand{\hsol}[1]{\widehat{\mathit{s}}(#1)}
\newcommand{\hsB}{\hsol{B_n}}

\newcommand{\hsD}{\hsol{D_n}}

\newcommand{\Asets}{[n{-}1]}
\newcommand{\Bsets}{\overline{[n{-}1]}}
\newcommand{\Dsets}{[n{-}1]'}

\newcommand{\equal}[1]{
{\stackrel{{\textstyle #1}}{\ {\textstyle =}\ } }}
\newcommand{\onto}{\twoheadrightarrow}
\newcommand{\map}[1]{\xrightarrow{#1}}

\newcommand{\inc}{\hookrightarrow}

\begin{document}

\title[Peak and descent algebras]{The peak algebra\\ and\\ the descent algebras
of types B and D} \author{Marcelo Aguiar}
\address{Department of Mathematics\\
         Texas A\&M University\\
         College Station\\
         Texas  77843\\
         USA}
\email{maguiar@math.tamu.edu}
\urladdr{http://www.math.tamu.edu/$\sim$maguiar}

\author{Nantel Bergeron}
\address{Department of Mathematics and Statistics\\
        York University \\
        Toronto, Ontario M3J 1P3\\
        Canada}
\email{bergeron@mathstat.yorku.ca}
\urladdr{http://www.mathstat.yorku.ca/bergeron}

\author{Kathryn Nyman}
\address{Department of Mathematics\\
         Texas A\&M University\\
         College Station\\
         Texas  77843\\
         USA}
\email{nyman@math.tamu.edu}
\urladdr{http://www.math.tamu.edu/$\sim$kathryn.nyman}

\thanks{The first author would like to thank Swapneel Mahajan for sharing his insight on descent
algebras and for interesting conversations.\\ Research of the second author supported in part by
CRC, NSERC and PREA}
\date{March 3, 2003}
\date{March 3, 2003}

\begin{abstract}
We show the existence of a unital subalgebra $\pppn$ of the symmetric group algebra linearly
spanned by  sums of permutations with a common peak set, which we call the peak
algebra.
We show that $\pppn$ is the image of the descent algebra of type B under the
map to the descent algebra of type A which forgets the signs, and also the image of the
descent algebra of type D.
The algebra $\pppn$ contains a two-sided ideal $\pppnint$ which is defined in terms of
{\em interior} peaks.
This object was introduced in previous work by Nyman~\cite{Nym}; we find that
it is the image of certain ideals of the descent algebras of types B and D
 introduced in~\cite{BB92a} and~\cite{BV}. We derive an exact sequence of the form
 $0\to\pppnint\to\pppn\to\ppp{n-2}\to 0$. We obtain this and many
other properties of the peak algebra and its peak ideal by first establishing
analogous results for signed permutations and then forgetting the signs.
In particular, we construct two new commutative semisimple subalgebras of the descent
algebra (of dimensions $n$ and $\ipartn+1$) by grouping permutations according to their number of
peaks or interior peaks.
We  discuss the Hopf algebraic structures that exist on the direct sums of the
spaces $\pppn$ and $\pppnint$ over $n\geq 0$ and explain the connection with previous
work of Stembridge~\cite{Ste}; we also obtain new properties of his {\em descents-to-peaks}
map and construct a type B analog.
\end{abstract}
\maketitle

\begin{figure}[!h]
\parbox{360pt}{\tiny\tableofcontents}
\end{figure}

\section*{Introduction} \label{S:intro}
A descent of a permutation $w\in\calSn$ is a position $i$ for which
$w_i>w_{i+1}$, while a peak is a position $i$ for which $w_{i-1}<w_i>w_{i+1}$.
The first in-depth study of the combinatorics of peaks was
carried out by Stembridge~\cite{Ste}, who developed an analog of Stanley's theory of
poset partitions where the notion of descents (of linear extensions of posets)
is replaced by the notion of peaks.  More recent works uncovered connections
between peaks and such disparate topics as the generalized
Dehn-Sommerville equations~\cite{BMSV00, BHV, ABS} and the
Schubert calculus of isotropic flag manifolds~\cite{BH, BS, BMSV02}.
Additional interest in the study of peaks grew from Nyman's thesis~\cite{Nym-t,Nym},
which showed that summing permutations according to their peak sets leads to a non-unital subalgebra
of the group algebra of the symmetric group, much as Solomon's descent algebra is
constructed in terms of descent sets. Work that followed includes~\cite{Sch,Hsi} and the
present paper.

Our work developed from the observation
that peaks of ordinary permutations are closely related to descents of
signed permutations (permutations of type B). This led us to the discovery
of a new object: a unital subalgebra of the descent algebra of the symmetric group,
which we call the {\em peak algebra}. This algebra is linearly spanned by sums of permutations
with a common set of peaks. The object considered in~\cite{Nym,Sch} is the linear
span of the sums of permutations with a common set of {\em interior} peaks.
We find that this is in fact a two-sided ideal of the peak algebra.
 We warn the reader that in~\cite{Nym,Sch} the terminology
``peak algebra'' refers to our peak ideal, viewed as a non-unital algebra (and our peak
algebra is not considered in those works).

We obtain several properties of the peak algebra and the peak ideal.
Our approach consists of obtaining more general properties that hold at the level
of signed permutations and then specializing by means
of the canonical map that forgets the signs. This allows for simpler proofs
of known results for the peak ideal, as well as many new properties for both
the peak algebra and the peak ideal.

We next review the contents of the paper in more detail.

Let $\SolA$, $\SolB$ and $\SolD$ denote Solomon's descent algebras of the Coxeter
systems of types A, B and D (see
Section~\ref{S:descent} for the basic definitions). Let $\calSn$, $\calBn$ and $\calDn$ be
the corresponding Coxeter groups. $\SolA$ is thus the
usual descent algebra of the symmetric group $\calSn$.

In Section~\ref{S:morphisms} we introduce a commutative diagram of morphisms of algebras
\[\xymatrix{{\SolB}\ar[rr]^{\chi}\ar[dr]_{\varphi} & &{\SolD}\ar[ld]^{\psi}\\
& {\SolA} }\]
The map $\chi$ is due to Mahajan~\cite{Mah01,Mah-t}. These maps exist at the level of groups
(for instance, $\varphi(w)$ is the permutation
obtained by forgetting the signs of the signed permutation $w$), but it is not obvious that
they restrict to descent algebras. We prove these facts in
Propositions~\ref{P:BDmap},~\ref{P:Btopeaks-D} and~\ref{P:Dtopeaks}. Central to this is
the analysis of the transformation of descents of signed permutations under the map $\varphi$, which
we carry out in Section~\ref{S:peaks}.
The principle that arises is: peaks are shadows of descents of type B and are not
straightforwardly related to descents of type A.

We will show that the images of $\varphi$ and $\psi$ coincide. This common image is
a certain subalgebra of $\SolA$ whose study is
the main goal of this paper. One of our main results (Theorem~\ref{T:peakalgebra})
states that this subalgebra is linearly
spanned by the elements
\begin{equation*}
P_F:=\sum_{\Peak(w)=F}w\,,
\end{equation*}
where $\Peak(w)$ denotes the set of peaks of a permutation $w$.
We refer to this subalgebra as the peak algebra and denote it by $\pppn$.
It is important to note that we allow $1$ as a possible peak (by making the convention
that $w_0=0$).

Previous work of Nyman~\cite{Nym-t,Nym} had shown the existence of a similar (but nonunital)
subalgebra of $\SolA$---the difference stemming solely from the fact that peaks at $1$ were not allowed
in that work. We obtain a stronger result in this work: we show that Nyman's algebra is in fact a two-sided ideal of $\pppn$
(Theorem~\ref{T:peakideal}).
This ideal, that we denote by $\pppnint$, is the image under $\varphi$
of the two-sided ideal $I^0_n$ of
$\SolB$ introduced in~\cite{BB92a}, and also under $\psi$ of a similar ideal of $\SolD$
introduced in~\cite{BV}. We obtain these results in
Theorems~\ref{T:Bideals}  and~\ref{T:Dexactsequence}.

We describe an exact sequence of the form
\[0\to\pppnint\to\pppn\to\ppp{n-2}\to 0\]
and relate it to the exact sequences of~\cite{BB92a} and~\cite{BV} involving the
descent algebras of types B and D (Corollary~\ref{C:Bexactsequence}
and Theorem~\ref{T:Dexactsequence}). This is an algebraic version of the Fibonacci recursion
$f_n=f_{n-1}+f_{n-2}$.

It is well-known that grouping permutations according to the number of descents leads to a
commutative semisimple subalgebra of
$\SolA$~\cite{Lod89}. In Section~\ref{S:commutative} we derive an analogous result for the linear span of the
sums of permutations with a given number of peaks. Once again, it is possible to easily derive
this fact from results for type B, which are known from~\cite{BB92b,Mah01}.
We find (Theorem~\ref{T:wp})
in fact two commutative semisimple subalgebras
$\wp_n$ and $\hwp_n$ of $\pppn$ and a two-sided ideal $\wpint_n$ of $\hwp_n$ such that
\[\hwp_n=\wp_n+\wpint_n\,.\]
In parallel with the situation for the peak algebra and its peak ideal,
 $\wp_n$ and $\wpint_n$ are obtained by grouping permutations according to the number of peaks and
  interior peaks, respectively.
Thus, we see that $\wpint_n$ coincides with the non-unital commutative subalgebra studied by
Schocker~\cite{Sch}.  Objects closely related to $\wp_n$ and $\wpint_n$ were first described
by Doyle and Rockmore~\cite{DR}. See Remarks~\ref{S:otherwork} for more details about these
connections.

In Corollary~\ref{C:sBexactsequence}, we describe an exact sequence relating $\wp_n$ and
$\wp_{n-2}$ and an analogous sequence at the level of type B. The corresponding objects
at the level of type D are discussed in Section~\ref{S:comm-D}.

Mantaci and Reutenaer have defined a certain subalgebra of the group algebra of $\calBn$
which strictly contains the descent algebra~\cite{ManReu}. We denote it by $\OmeB$. We recall its
definition in Section~\ref{S:ManReu}, and provide several new results that link $\OmeB$ to
$\SolA$ by means of $\varphi$. The correspondence between types B and A is then summarized
by the diagram
\[\xymatrix@C=0pc{ {I_n^{0}}\ar@{->>}[d]_{\varphi} &\subseteq &
{\SolB}\ar@{->>}[d]_{\varphi} &\subseteq& {\OmeB}\ar@{->>}[d]^{\varphi} &\subseteq&
 {\Q\calBn}\ar@{->>}[d]^{\varphi}\\
 {\pppnint} &\subseteq& {\pppn}& \subseteq & {\SolA} &\subseteq& {\Q\calSn} }\]
We show that $I^0_n$ is a right ideal
of $\OmeB$ and that it is in fact principal (Corollary~\ref{C:rightOmeB}).
This allows us to recover a result of Shocker:
that $\pppnint$ is a principal right ideal of $\SolA$ (Corollary~\ref{C:rightSolA}).
In addition, we obtain a new result
that states that $\pppnint$ is principal as a right ideal of $\pppn$ (this
result neither implies Schocker's result nor is it implied by it), as well as the corresponding
statement for type B (Proposition~\ref{P:rightpeak}).

An important map $\Theta:\SolA\to\pppnint$, which we call the {\em descents-to-peaks transform},
is discussed in Section~\ref{S:thetamaps}. We construct
a type B analog $\ThetaB:\OmeB\to I^0_n$, from which the basic properties of $\Theta$ may
be easily derived. The descent-to-peaks transform is a special case of a map considered
by Krob, Leclerc and Thibon~\cite{KLT} and is dual to a map considered by Stembridge~\cite{Ste}. See
Remarks~\ref{R:thetadef} for the precise details.

In Section~\ref{S:external} we consider the direct sum over all $n\geq 0$ of the group, descent and
peak algebras discussed in previous sections. This leads to a diagram
\[\xymatrix@C=0pc{ {I^{0}}\ar@{->>}[d]_{\varphi} &\subseteq &
{\Sol{B}}\ar@{->>}[d]_{\varphi} &\subseteq& {\Ome{B}}\ar@{->>}[d]^{\varphi} &\subseteq&
 {\QB}\ar@{->>}[d]^{\varphi}\\
 {\pppint{}} &\subseteq& {\ppp{}}& \subseteq & {\Sol{A}} &\subseteq& {\QS} }\]
It is known from work of Malvenuto and Reutenauer~\cite{MalReu} that the space $\QS$  may be endowed with a new
product (called the external product) and a coproduct that turn it into a graded Hopf algebra.
Moreover, this structure restricts to $\Sol{A}$, which results in the Hopf algebra
of non-commutative symmetric functions. Using results from~\cite{AM},
we discuss type B analogs of these constructions
and show that all objects in the above diagram are Hopf algebras, except for $\Sol{B}$,
which is an $I^0$-module coalgebra, and $\ppp{}$, which is a $\pppint{}$-module coalgebra.
We also discuss the behavior with respect to the external structure of all maps
from previous sections.

In Section~\ref{S:duality} we clarify the connection between Stembridge's Hopf algebra of
peaks and the objects discussed in this paper. The Hopf algebra $\pppint{}$ and
the descents-to-peaks transform are dual to the objects considered by Stembridge in~\cite{Ste}.

In Section~\ref{S:words} we recall the canonical action of permutations on words, and
provide an explicit description for the action of the generators of the
principal ideals $I^0_n$ and $\pppnint$ in terms of {\em symmetrizers} and {\em Jordan brackets}
(Propositions~\ref{P:action-B} and~\ref{P:action-P}).
This complements a result of Krob, Leclerc and Thibon~\cite{KLT}.

Most of the results in this paper were presented by one of us (Nyman) at the
meeting of the American Mathematical Society in Montr\'eal in May, 2002.

\section{Descent algebras of Coxeter systems}\label{S:descent}

Let $(W,S)$ be a finite Coxeter system~\cite{Hum}. That is, $W$ is a finite group generated
by the set $S$ subject to the relations
\[(st)^{m_{st}}=1 \text{ for all }s,t\in S\,,\]
where the $m_{st}$ are positive integers and $m_{ss}=1$ for all $s\in S$.

Given $w\in W$, its descent set is
\[\Des(w):=\{s\in S\ \mid\ \ell(ws)<\ell(w)\}\,,\]
where $\ell(w)$ denotes the length of a minimal expression for $w$ as a product
of elements of $S$.

Let $\Q W$ denote the group algebra of $W$. The subspace spanned by the
elements
\begin{equation}\label{E:defDbasis}
Y_J:=\sum_{\Des(w)=J}w
\end{equation}
is closed under the product of $\Q W$~\cite{Sol}. It is the Solomon's descent algebra
of $(W,S)$.   The unit element is $Y_\emptyset=1$.

The set $\{Y_J\}_{J\subseteq S}$ is a linear basis of the descent algebra.
A second linear basis is defined by
\begin{equation}\label{E:defBbasis}
X_J:=\sum_{I\subseteq J}Y_I=\sum_{\Des(w)\subseteq J}w\,.
\end{equation}
Note that
\begin{equation}\label{E:DfromB}
Y_J=\sum_{I\subseteq J}(-1)^{\#J-\#I}X_I
\end{equation}
so the set $\{X_J\}_{J\subseteq S}$ indeed forms a basis. This basis has proved useful
in describing the structure of descent algebras~\cite{Sol, GRem, GReu, BB92a, BV} and
will be useful in our work as well.

{\bf Warning}: The notations $Y_J$ and $X_J$ do not make reference to the Coxeter system,
which will have to be understood from the context. This is particularly
relevant in our work in which we deal with descent algebras of
various Coxeter systems at the same time.

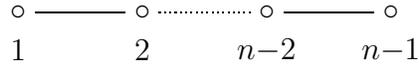
\begin{figure}[htb]
\[\xymatrix@C+10pt{
 {\circ}\ar@{-}[r] \save[]-<0cm,0.5cm>*{1}\restore
 &{\circ}\ar@{.}[r] \save[]-<0cm,0.5cm>*{2}\restore
 &{\circ}\ar@{-}[r] \save[]-<0cm,0.5cm>*{n{-}2}\restore
&{\circ}\save[]-<0cm,0.5cm>*{n{-}1}\restore
 }\]
\caption{The Coxeter graph $A_{n-1}$\label{F:Agraph}}
\end{figure}

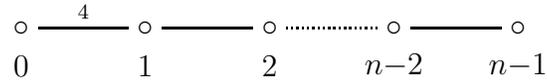
\begin{figure}[htb]
 \[\xymatrix@C+10pt{
 {\circ}\ar@{-}[r]^{4} \save[]-<0cm,0.5cm>*{0}\restore
 &{\circ}\ar@{-}[r] \save[]-<0cm,0.5cm>*{1}\restore
 & {\circ}\ar@{.}[r] \save[]-<0cm,0.5cm>*{2}\restore
 &{\circ}\ar@{-}[r] \save[]-<0cm,0.5cm>*{n{-}2}\restore
&{\circ}\save[]-<0cm,0.5cm>*{n{-}1}\restore
 }\]
  \caption{The Coxeter graph $B_{n}$\label{F:Bgraph}}
\end{figure}

\begin{figure}[htb]
 \[\xymatrix@C+10pt{
 {\circ}\ar@{-}[rd] \save[]-<0cm,0.5cm>*{1}\restore\\
  &{\circ}\ar@{-}[r]\save[]-<0cm,0.5cm>*{2}\restore
  &{\circ}\ar@{.}[r] \save[]-<0cm,0.5cm>*{3}\restore
  &{\circ}\ar@{-}[r] \save[]-<0cm,0.5cm>*{n{-}2}\restore
  &{\circ} \save[]-<0cm,0.5cm>*{n{-}1}\restore\\
{\circ}\ar@{-}[ru] \save[]-<0cm,0.5cm>*{1'}\restore
 }\]
  \caption{The Coxeter graph $D_{n}$\label{F:Dgraph}}
\end{figure}
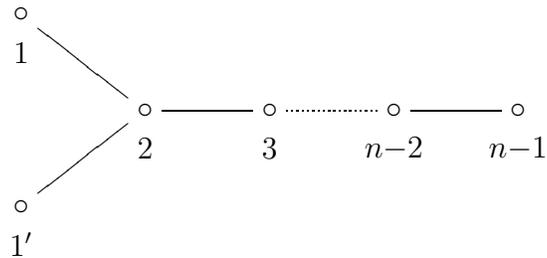

From now on, we assume that the Coxeter system $(W,S)$ is associated to a
Coxeter graph $G$ : $S$ is the set of vertices and $m_{st}$ is the label of the edge
joining $s$ and $t$ ($2$ if there is no such edge and $3$ is there is an edge with no
label, see Figures~\ref{F:Agraph}, \ref{F:Bgraph} and \ref{F:Dgraph}). We denote the descent algebra
of this Coxeter system by $\Sol{G}$.

For the Coxeter graph $A_{n-1}$, the Coxeter group is the symmetric group $\calSn$
and the set of generators consists of the elementary transpositions $s_i=(i,i+1)$
for $i=1,\ldots,n-1$. We identify this set with the set $\Asets:=\{1,2,\ldots,n{-}1\}$
via $s_i\leftrightarrow i$.
A permutation $w$ is represented by a sequence $w=w_1\ldots w_n$ of distinct symbols from the
set $[n]$, where $w_i=w(i)$ are the values of the permutation.
A permutation $w\in \calSn$ has a descent at $i\in\Asets$ if $w_i>w_{i+1}$.

 For the Coxeter graph $B_n$, the Coxeter group is the group of signed permutations
$\calSn\ltimes\Z_2^n$, which we denote by $\calBn$. A signed permutation $w$
is represented by a sequence $w=w_1\ldots w_n$ of symbols from $1$ to $n$ (the underlying permutation from
$\calSn$), some of which may be barred (according to the sign from $\Z_2^n$). We order
these symbols by
\begin{equation} \label{E:order}
\ldots<\bar{2}<\bar{1}<1<2<\ldots\,.
\end{equation}
The set of generators consists of the elementary transpositions $s_i$ as before
(with no signs) plus the signed permutation
\[s_0:=\bar{1}2\ldots n\,.\]
We identify this set with the set $\Bsets:=\{0\}\cup\Asets$.
A signed permutation $w\in \calBn$ has a descent at $i\in\Bsets$ if
$w_i>w_{i+1}$, with respect to the order~\eqref{E:order}, where we agree that $w_0=0$.

For the Coxeter graph $D_n$, the Coxeter group is the subgroup of $\calBn$ consisting of
those signed permutations with an even number of signs (bars). We denote it by
 $\calDn$.

The set of generators consists of the elementary transpositions $s_i$ as before,
from $i=1$ to $n{-}1$, plus the signed permutation
\[s_{1'}:=\bar{2}\bar{1}3\ldots n \,.\]
We identify this set with the set $\Dsets:=\{1'\}\cup\Asets$.

It is convenient to represent a signed permutation
 $w=w_1\ldots w_n\in\calDn$ by a fork-shaped  sequence
 \[\xymatrix@C=-5pt@R=-5pt{
 {\!w_1}\ar@{-}[rd] \\
  &{w_2} \save[]-<2cm,0cm>*{w\ \ =}\restore
  &{w_3} &{\cdots}  &{w_{n{-}1}} &{w_{n}} \\
 {w_{1'}}
 }\]
of symbols from $1$ to $n$, possibly barred, with the convention that
$w_{1'}=\overline{w_1}$.
It then turns out that an element
$w\in \calDn$ has a descent at $i\in \Dsets$ if $w_i>w_{i+1}$, with respect to
the order $(*)$, and where we understand that $1'+1=2$. For instance,
 \[\xymatrix@C=-3pt@R=-5pt{
 {\!\bar{2}} \\
  &{\bar{1}} \save[]-<1.2cm,0cm>*{s_{1'}\ \ =\ \ \ \ }\restore
  &{3}&\cdots &{n{-}1} &{n} \\
 {2}
 }\]
has a descent at $1'$ because $2>\bar{1}$, but not at $1$, since $\bar{2}<\bar{1}$.

\section{Morphisms between descent algebras}\label{S:morphisms}

There is a canonical morphism $\varphi:\calBn\onto \calSn$ from the group of signed permutations onto
that of ordinary permutations, obtained by simply forgetting the signs (bars).
Let $\psi:\calDn\onto \calSn$ be its restriction. Consider the linear extensions of
these maps to the group algebras. We will show in Section~\ref{S:peaks} that these maps restrict to
the corresponding descent algebras, yielding morphisms of algebras
\[\varphi:\SolB\to\SolA \text{ \ and \ }\psi:\SolD\to\SolA\,.\]
It turns out that there is also a morphism of algebras
\[\chi:\SolB\to\SolD\]
such that $\varphi=\psi\chi$. In this section we concentrate on this map.

The map $\chi$ is due to Mahajan.
It is {\em not} the restriction of a morphism of groups $\calBn\to \calDn$, and
its construction is better understood from the point of view of hyperplane
arrangements, as explained in~\cite[Section 6.1]{Mah01}. We will add a group-theoretic counterpart
to Mahajan's geometric construction  that makes the above
commutativity evident (but does not explain why $\chi$ is a morphism of algebras).

We begin by defining a map $\chi:\calBn\onto \calDn$ as
\[\chi(w)=\begin{cases} w &\text {if }w\in \calDn\,,\\
ws_0 &\text{ if }w\notin \calDn\,. \end{cases}\]
In other words, $\chi$ changes the sign of the first entry if the number of signs is odd.
Note that $\chi$ is not a morphism of groups. Nevertheless, we have the following

\bp\label{P:BDmap} The map $\chi$ restricts to descent algebras and, at that level,
it is a morphism of algebras $\chi:\SolB\to\SolD$. Moreover,
it is explicitly given as follows:
for  $J\subseteq\{2,\ldots,n{-}1\}$,
\begin{align}
Y_J & \mapsto Y_J\notag\\
Y_{\{1\}\cup J} & \mapsto Y_{\{1'\}\cup J}+ Y_{\{1\}\cup J}+Y_{\{1',1\}\cup J} \notag\\
Y_{\{0\}\cup J} & \mapsto Y_J+Y_{\{1'\}\cup J}+ Y_{\{1\}\cup J} \label{E:BDmap1}\\
Y_{\{0,1\}\cup J} &  \mapsto  Y_{\{1',1\}\cup J}\notag\\
\intertext{or equivalently,}
X_J & \mapsto X_J\notag\\
X_{\{1\}\cup J} & \mapsto X_{\{1',1\}\cup J} \notag\\
X_{\{0\}\cup J} & \mapsto X_{\{1'\}\cup J}+ X_{\{1\}\cup J}\label{E:BDmap2}\\
X_{\{0,1\}\cup J} &  \mapsto 2\cdot X_{\{1',1\}\cup J}\notag
\end{align}
\ep
\bpf Consider $Y_{\{1\}\cup J}\in\SolB$. Take $w\in\calBn$ such that
$\Des(w)=\{1\}\cup J$.  There are exactly three possibilities for $\Des(\chi(w))$:
\begin{quote}
if $-w_1>w_2$ then $\Des(\chi(w))= {\{1',1\}\cup J}$,\\
if $-w_1<w_2$ and $w\in\calDn$ then $\Des(\chi(w))= {\{1\}\cup J}$, and\\
if $-w_1<w_2$ and $w\notin\calDn$ then $\Des(\chi(w))= {\{1'\}\cup J}$.
\end{quote}
Conversely, given $v\in\calDn$ such that  $\Des(v)={\{1\}\cup J}$, we must have $v_1>0$;
so letting $w=v\in\calBn$ we obtain $\Des(w)=\{1\}\cup J$ and $\chi(w)=v$.
If, instead, $\Des(v)={\{1'\}\cup J}$, then we must have $v_1<0$, so letting
$w=vs_0\in\calBn$ we still obtain $\Des(w)=\{1\}\cup J$ and $\chi(w)=v$.
Finally, if $\Des(v)={\{1',1\}\cup J}$, then we let $w=v\in\calBn$ if $v_1>0$
and $w=vs_0\in\calBn$ if $v_1<0$, to obtain the same conclusion.

Therefore,
\[\chi(Y_{\{1\}\cup J})=Y_{\{1'\}\cup J}+ Y_{\{1\}\cup J}+Y_{\{1',1\}\cup J}\,.\]
The remaining expressions in \eqref{E:BDmap1} may be similarly deduced.
The expressions in \eqref{E:BDmap2} follow easily from \eqref{E:defBbasis}.

Comparing~\eqref{E:BDmap2} to~\cite[Section 6.1]{Mah01}
we see that $\chi$ coincides with Mahajan's map.
His construction guarantees that $\chi:\SolB\to\SolD$ is a morphism of algebras.
\epf

\begin{rem} It is possible to give a direct proof of the fact that $\chi$ is a morphism of algebras
by making use of the multiplication rules for the $X$-bases of
$\SolB$ and $\SolD$ given in~\cite[Theorem 1]{BB92a} and~\cite[Theorem 1]{BV}.
\end{rem}

\bc\label{C:BDA} There is a commutative diagram of morphisms of algebras
\begin{equation}\label{E:BDA}
\xymatrix{{\SolB}\ar[rr]^{\chi}\ar[dr]_{\varphi} & &{\SolD}\ar[ld]^{\psi}\\
& {\SolA} }
\end{equation}
\ec
\bpf The commutativity $\varphi=\psi\chi$ already holds at the level of groups, clearly.
\epf

The descent algebras $\SolB$ and $\SolD$ have the same dimension $2^n$, but the
map $\chi$ is not an isomorphism: its image has codimension $2^{n-2}$.
In fact, the image of $\chi$ may be easily described. Consider the
decomposition $\calDn=\calDn^{(1)}\coprod\calDn^{(2)}\coprod\calDn^{(3)}$
where the $\calDn^{(i)}$ are respectively defined  by the conditions
\[\abs{w_1}<w_2,\ \ \abs{w_1}>\abs{w_2} \text{ \ and \ }\abs{w_1}<-w_2\,.\]
In the group algebra of $\calDn$, define elements
\[Y_J^{(i)}=\sum\{w\in\calDn^{(i)}\ \mid\ \Des(w)\cap\{2,\ldots,n{-1}\}=J\},\]
for each $J\subseteq\{2,\ldots,n{-}1\} \text{ and }i=1,2,3$.

\bc\label{C:imchi}  The image of $\chi:\SolB\to\SolD$ is linearly spanned by
the elements $Y_J^{(i)}$.
In particular, these elements span a subalgebra of $\SolD$ of dimension
$3\cdot 2^{n-2}$ (and codimension $2^{n-2}$).
\ec
\bpf It suffices to observe that
\[
 Y_J^{(1)}=Y_J,\ \ 
  Y_J^{(2)}=Y_{\{1'\}\cup J}+ Y_{\{1\}\cup J} \text{ \ and \ }
  Y_J^{(3)}=Y_{\{1',1\}\cup J}\]
and to make use of \eqref{E:BDmap1}.
\epf

The maps $\varphi$ and $\psi$ are harder to describe in terms of the $Y$ or $X$ bases.
An explicit description is provided below (Propositions~\ref{P:Btopeaks-D}, ~\ref{P:Btopeaks-B}
and~\ref{P:Dtopeaks}). It turns out that these maps have the same image, even though $\chi$ is
not surjective. The main goal of this paper is to describe this image, a certain
subalgebra of $\SolA$, in explicit combinatorial terms.

\begin{rem} A similar commutative diagram to~\eqref{E:BDA} is given in~\cite[Section 6]{Mah01}.
In fact, the map $\chi$ is the same for both, but the maps $\SolB\to\SolA$ and $\SolD\to\SolA$
considered in~\cite{Mah01} are different from ours.
\end{rem}

\section{Peaks of permutations and descents of signed permutations} \label{S:peaks}

\bd\label{D:peaks}
Let $w\in\calSn$. The set of {\em peaks} of $w$ is
\[\Peak(w)=\{i\in\Asets\ \mid\ w_{i-1}<w_i>w_{i+1}\}\,,\]
where we agree that $w_0=0$.
\ed
For instance, $\Peak(42153)=\{1,4\}$.
Note that $1$ is a peak of $w$ if and only if it is a descent of $w$.
Thus, our notion of peaks differs (slightly) from that of~\cite{Ste,Nym,Sch}.
This turns out to be an important distinction, as will be made clear throughout this
work. We will deal with both notions of peaks starting in Section~\ref{S:ideals}.

Let $\calQn$ denote the collection of all subsets of $\Asets$, and
\[\calFn=\{F\in\calQn\ \mid\ \text{if $i\in F$ then $i+1\notin F$}\}.\]
Throughout this work, by the $n$-th Fibonacci number we understand
\begin{equation}\label{E:fibonacci}
f_n:=\#\calFn\,.
\end{equation}
Thus, $f_0=1$, $f_1=1$, $f_2=2$ and $f_n=f_{n-1}+f_{n-2}$.

Clearly, the peak set of any permutation $w\in \calSn$ belongs to $\calFn$.
It is easy to see that any element of $\calFn$ is realized in this way.

Given $J\in\calQn$, consider the set
\[\Lambda(J):=\{i\in J\ \mid\ i-1\notin J\}.\]
One verifies immediately that $\Lambda(J)\in\calFn$ and, moreover, that
the diagram
\begin{equation}\label{E:peak-descent}
\smallxymatrix{ & {\calSn}\ar[dl]_{\Des}\ar[dr]^{\Peak} & \\
{\calQn}\ar[rr]_{\Lambda}& & {\calFn} }
\end{equation}
commutes.

For each $F\in\calFn$, define an element of the group algebra of $\calSn$ by
\begin{equation}\label{E:defpeakone}
P_F:=\sum_{\Peak(w)=F}w\,.
\end{equation}
It follows from~\eqref{E:defDbasis} and~\eqref{E:peak-descent} that
\begin{equation}\label{E:defpeaktwo}
P_F=\sum_{\Lambda(J)=F}Y_J\,;
\end{equation}
in particular, each $P_F$ belongs to the descent algebra $\SolA$.

 It turns out that $\varphi:\SolB\to\SolA$ maps a basis element $Y_J$ to a certain sum of
elements $P_F$. In order to derive an explicit expression for $\varphi(Y_J)$, we need
to introduce some notation and make some basic observations.

First, given a permutation $u\in\calSn$ and a position $i\in[n]$,
we refer to one of the four diagrams
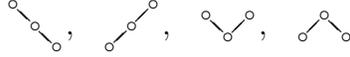
\[\descent,\ \  \ascent,\ \  \valley,\ \  \peak\]
as the {\em local shape of $u$ at $i$}, according to whether $u_{i-1}>u_i>u_{i+1}$,
$u_{i-1}<u_i<u_{i+1}$, $u_{i-1}>u_i<u_{i+1}$ or $u_{i-1}<u_i>u_{i+1}$, respectively.
Here, we make the conventions that $u_0=0$ and $u_{n+1}=+\infty$.
A position $i$ is a peak of $u$ if the local shape at $i$ is $\peak$; we also say
that $i$ is a {\em valley} if the local shape is $\valley$.
Note that valleys may only occur in positions $2,\ldots,n$, while peaks only in positions
$1,\ldots,n{-}1$. Note also that the total number of peaks equals the total number of valleys,
since in listing them from left to right, peaks and valleys alternate, starting
always with a peak and ending with a valley (except for the case of the identity permutation,
which has no peaks or valleys). For instance, the four local shapes for  $u=2413$ are, from left
to right,
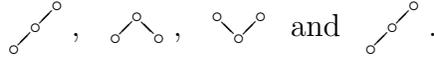
\[\ascent\,,\ \  \peak\,,\ \ \valley \text{ \ and \ }\ascent \,.\]

Similarly, for a signed permutation
$w\in\calBn$ and $i\in[n]$, we speak of its local shape at $i$, which is one of
the previous diagrams decorated by the sign of the entry $w_i$, with the same conventions.
For instance, the four local shapes of $w=2\bar{4}13$ are (in order)
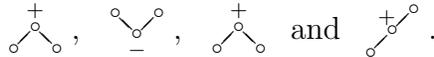
\[\speak{+}\,,\ \  \svalley{-}\,,\ \ \speak{+} \text{ \ and \ }\sascent{+} \,.\]

Clearly, the collection of all local shapes of $w\in\calBn$ determine its descent set,
and conversely. More precisely, given a subset $J\subseteq\Bsets$, we have that $\Des(w)=J$
if and only if for each $i\in[n]$, the local shape of $w$ at $i$ is
as follows
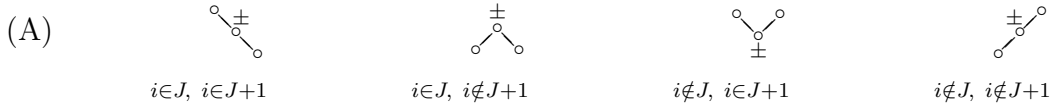
\begin{align*}\tag{A}
\sdescent{\pm} & &\speak{\pm}  & &\svalley{\pm} & & \sascent{\pm}\\
\scriptstyle{ i\in J,\ i\in J+1} & &\scriptstyle{i\in J,\ i\notin J+1} & &
\scriptstyle{i\notin J,\ i\in J+1} & &\scriptstyle{i\notin J,\ i\notin J+1}
\end{align*}
where $J+1$ denote the shifted set $J+1=\{i+1\ \mid\ i\in J\}$.

Now suppose that a permutation $u\in\calSn$ is given and we look for all $2^n$ signed permutations
$w\in\calBn$ such that $\varphi(w)=u$. Given a local shape of $u$ at $i\in[n]$, the possible
local shapes of $w$ at $i$ are as follows
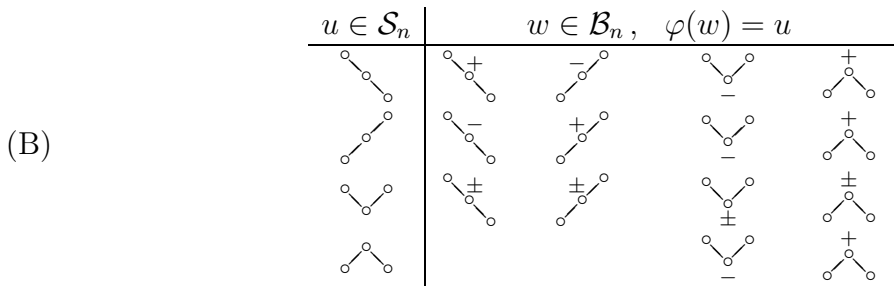
\begin{equation}\tag{B}
\begin{tabular}{c|cccc}
$u\in\calSn$ &  & $w\in\calBn$\,, & $\varphi(w)=u$& \\
\hline
$\descent$ & $\sdescent{+}$  & $\sascent{-}$ & $\svalley{-}$ & $\speak{+}$\\
$\ascent$ & $\sdescent{-}$  & $\sascent{+}$ & $\svalley{-}$ & $\speak{+}$\\
$\valley$ & $\sdescent{\pm}$  & $\sascent{\pm}$ & $\svalley{\pm}$ & $\speak{\pm}$\\
$\peak$ &   &  & $\svalley{-}$ & $\speak{+}$
\end{tabular}
\end{equation}

Let $I\vartriangle J=(I\cup J)\setminus(I\cap J)$ be the symmetric set difference.

\bp\label{P:Btopeaks-D}
The map $\varphi:\SolB\to\SolA$ is explicitly given by
\begin{equation}
Y_J \mapsto \sumsub{F\in\calFn\\F\subseteq J\vartriangle (J+1)}2^{\# F}\cdot P_F\,,
\label{E:Btopeaks-D}
\end{equation}
for any $J\subseteq\Bsets$.
\ep
\bpf We need to show that for each $u\in\calSn$,
\begin{equation} \tag{$\ast$}
\#\{w\in\calBn\ \mid\ \Des(w)=J \text{ and
}\varphi(w)=u\}=\begin{cases} 2^{\#\Peak(u)} & \text{ if }\Peak(u)\subseteq J\vartriangle(J+1)\,,\\
0 & \text{ if not.}
\end{cases}
\end{equation}
Note that this set is defined precisely by conditions (A) and (B) above.

Suppose there is at least one  $w$ in the
above set and take $i\in\Peak(u)$. By (B), the local shape of $w$ at $i$ is either
$\svalley{-}$ or $\speak{+}$. Then, by (A),  $i\in J\vartriangle(J+1)$. This proves
the second half of ($\ast$).

Assume now that $\Peak(u)\subseteq J\vartriangle(J+1)$.
A signed permutation $w$ in the
above set is completely determined once the signs of its entries $\abs{w_i}=u_i$ are chosen. We will show
that the signs at those positions $i$ that are valleys of $u$ can be freely chosen,
while the signs at the remaining positions are uniquely determined.
This proves the remaining half of ($\ast$), because as explained above, the number
of valleys equals the number of peaks, so the number of such choices is $2^{\#\Peak(u)}$.

Given a position $i$, it will fall into one of the four alternatives of (A) (according to $J$)
and one of the four alternatives of (B) (according to $u$).
Consider first the case of a position $i$ that is not a valley of $u$. In this case,
whatever the alternative in (A) is, there is exactly one possible local shape
for $w$ that satisfies (B). In particular,
the sign of $w_i$ is uniquely determined for all those positions $i$.
On the other hand, if $i$ is a valley of $u$, then whatever the alternative in (A) is,
there are exactly two possible local shapes for $w$ that satisfy (B), each with the same
underlying shape but with different signs.
Thus, the choices of signs that lead to a signed permutation
$w$ that satisfies (A) and (B) is the same as the choices of signs on the set of
valleys of $u$, as claimed.
\epf

We can  derive a similar expression for $\varphi$ on the $X$-basis of $\SolB$.
\bp\label{P:Btopeaks-B}
The map $\varphi:\SolB\to\SolA$ is explicitly given by
\begin{equation}
X_J \mapsto 2^{\# J}\cdot\!\!\!\!\sumsub{F\in\calFn\\ F\subseteq J\cup(J+1)} P_F\label{E:Btopeaks-B}\\
\end{equation}
for any $J\subseteq\Bsets$.
\ep
\bpf Since $X_J=\sum_{K\subseteq J}Y_K$, equation~\eqref{E:Btopeaks-B} will follow from~\eqref{E:Btopeaks-D},
once we show that for any $F\in\calFn$ and $J\subseteq\Bsets$ we have
\[\#\{K\subseteq\Bsets\ \mid\ F\subseteq K\vartriangle(K+1),\ K\subseteq J\}=
\begin{cases} 2^{\#J-\#F} & \text{ if }F\subseteq J\cup(J+1),\\
0 & \text{ if not.} \end{cases}\]
If there is at least one subset $K$ in the above family, then
$F\subseteq K\cup(K+1)\subseteq J\cup(J+1)$. This proves the second case of the equality.

Suppose now that $F\subseteq J\cup(J+1)$. We decompose an arbitrary subset $K$ of $J$ into five pieces:
\[K^0:=K\cap\Bigl(J\setminus\bigl(F\cup(F-1)\bigr)\Bigr)\,,\ \
K^+:=K\cap F\ \ \text{and}\ \ K^-:=K\cap(F-1)\,,\]
\[K^+_1:=K\cap F\cap(J+1)\,,\ \ K^+_2:=K\cap F\cap(J+1)^c\,,\]
\[K^-_1:=K\cap (F-1)\cap(J-1)\,,\ \ K^-_2:=K\cap (F-1)\cap(J-1)^c\,.\]
Thus, $K=K^0\sqcup K^+_1\sqcup K^+_2\sqcup K^-_1\sqcup K^-_2$.
One may check that, under the present assumptions ($F\in\calFn$ and $F\subseteq J\cup(J+1)$), a subset $K$ of $J$ satisfies
$F\subseteq K\vartriangle(K+1)$  if and only if the following conditions hold:
\[K^+_2=J\cap F\cap(J+1)^c\,,\ \ K^-_2=J\cap (F-1)\cap(J-1)^c \ \ \text{ and }\ \
K^+_1\sqcup(K^-_1+1)=J\cap F\cap(J+1)\,,\]
while $K^0$ may be any subset of $J\setminus\bigl(F\cup(F-1)\bigr)$ and $K^+_1$ any subset
of $J\cap F\cap(J+1)$.  The claim now follows from the fact that
\[\#J\setminus\bigl(F\cup(F-1)\bigr)+\#J\cap F\cap(J+1)=\#J-\#F\]
(this also makes use of the assumptions $F\in\calFn$ and $F\subseteq J\cup(J+1)$).
\epf
\begin{rem} Formulas~\eqref{E:Btopeaks-D} and~\eqref{E:Btopeaks-B} are very similar to certain
expressions for a map introduced by Stembridge~\cite[Propositions 2.2, 3.5]{Ste}.
This is in fact a consequence of a very non-trivial fact: that the dual of Stembridge's
map admits a {\em type B analog} and $\varphi$ commutes with these maps.
This connection will be clarified in Sections~\ref{S:thetamaps} and~\ref{S:duality}.
\end{rem}

\begin{rem}\label{R:kerphi} It is interesting to note that for any $J\subseteq\Bsets$,
\[\varphi(Y_J)=\varphi(Y_{\Bsets\setminus J})\,.\]
This can be deduced from the explicit formula~\eqref{E:Btopeaks-D},
but it can also be understood as follows. Consider the map $\sigma:\calBn\to \calBn$
which reverses the sign of all the entries of a signed permutation ($\sigma$ is
both left and right multiplication by $\bar{1}\bar{2}\ldots\bar{n}$).  Clearly,
\[\Des\bigl(\sigma(w)\bigr)=\Bsets\setminus\Des(w)\,,
\text{ and therefore }\sigma(Y_J)=Y_{\Bsets\setminus J}\,.\]
On the other hand, since $\varphi$ simply forgets the signs, we have
\[\varphi(w)=\varphi\bigl(\sigma(w)\bigr)\,, \text{ and therefore }
\varphi(Y_J)=\varphi\bigl(\sigma(Y_J)\bigr)=
\varphi(Y_{\Bsets\setminus J})\,.\]
\end{rem}

\medskip
The map $\psi:\SolD\to\SolA$ also sends basis elements $Y_J$ and $X_J$ to sums of
elements $P_F$. The next result provides explicit expressions.

We make use of the convention $1'+1=2$ to give meaning to the shifted set $J+1=\{i+1\ \mid\ i\in J\}$,
where $J$ is any subset of $\Dsets$.

\bp\label{P:Dtopeaks}  The map $\psi:\SolD\to\SolA$ is explicitly given by
\begin{align}\label{E:Dtopeaks-D}
Y_J & \mapsto \sumsub{F\in\calFn\\F\subseteq J\vartriangle (J+1)}2^{\# F}\cdot P_F \notag\\
Y_{\{1\}\cup J} \text{ and }Y_{\{1'\}\cup J} & \mapsto
\sumsub{\{1\}\cup F\in\calFn\\F\subseteq J\vartriangle (J+1)}2^{\# F}\cdot
P_{\{1\}\cup F}\\
Y_{\{1',1\}\cup J} & \mapsto \sumsub{F\in\calFn\\F\subseteq J\vartriangle
\bigl(\{2\}\cup (J+1)\bigr)}2^{\# F}\cdot P_F \notag
\end{align}
and also
\begin{align}\label{E:Dtopeaks-B}
X_J & \mapsto 2^{\# J}\cdot\!\!\!\!\sumsub{F\in\calFn\\ F\subseteq J\cup(J+1)} P_F \notag\\
X_{\{1\}\cup J} \text{ and } X_{\{1'\}\cup J}& \mapsto 2^{\# J}\cdot\!\!\!\!\sumsub{F\in\calFn\\
F\subseteq J\cup(J+1)\cup\{1\}} P_F\\
X_{\{1',1\}\cup J} & \mapsto 2^{\# J+1}\cdot\!\!\!\!\sumsub{F\in\calFn\\
F\subseteq J\cup(J+1)\cup\{1,2\}} P_F \notag
\end{align}
for any $J\subseteq\{2,\ldots,n{-}1\}$.
\ep
\bpf It is possible to derive these expressions by a direct analysis of how the notion of descents
of type D transforms after forgetting the signs, as we did for the type B case in
Proposition~\ref{P:Btopeaks-D}. We find it much easier, however, to derive these from the analogous
expressions for the type B case, by means of the commutativity of diagram~\eqref{E:BDA}.

If $J\subseteq\{2,\ldots,n{-}1\}$, then by Proposition~\ref{P:BDmap}, the map $\chi:\SolB\to\SolD$
satisfies $\chi(Y_J)=Y_J$. Therefore, by Corollary~\ref{C:BDA} and Proposition~\ref{P:Btopeaks-D},
\[\psi(Y_J)=\varphi(Y_J)=
\sumsub{F\in\calFn\\F\subseteq J\vartriangle (J+1)}2^{\# F}\cdot P_F\,.\]
This is the first formula in~\eqref{E:Dtopeaks-D}. Similarly, the last formula in~\eqref{E:Dtopeaks-D}
follows from the facts that $\chi(Y_{\{0,1\}\cup J})=Y_{\{1',1\}\cup J}$ and
\[ \bigl(\{0,1\}\cup J\bigr)\vartriangle\Bigl(\bigl(\{0,1\}\cup J\bigr)+1\Bigr)=\{0\}\cup\Bigl(J\vartriangle
\bigl(\{2\}\cup (J+1)\bigr)\Bigr)\]
(note that $\{0\}$ does not affect the summation over subsets $F\in\calFn$).

The formulas for $Y_{\{1\}\cup J}$ and $Y_{\{1'\}\cup J}$ require a little extra work.
First, from Proposition~\ref{P:BDmap} we know that
\[\chi(Y_{\{0\}\cup J}-Y_J)= Y_{\{1'\}\cup J}+Y_{\{1\}\cup J}\,;\]
hence, by Corollary~\ref{C:BDA},
\[\psi(Y_{\{1'\}\cup J}+Y_{\{1\}\cup J})=\varphi(Y_{\{0\}\cup J}-Y_J)\,.\]
Now, since $J\subseteq \{2,\ldots,n-1\}$, $(\{0\}\cup J)\vartriangle\bigl((\{0\}\cup J)+1\bigr)=
\{0,1\}\cup\bigl(J\vartriangle(J+1)\bigr)$. Therefore, by Proposition~\ref{P:Btopeaks-D},
\begin{align*}
\varphi(Y_{\{0\}\cup J}-Y_J) &=
\sumsub{F\in\calFn\\F\subseteq \{1\}\cup\bigl(J\vartriangle (J+1)\bigr)}2^{\# F}\cdot
P_F-\sumsub{ F\in\calFn\\F\subseteq J\vartriangle (J+1)}2^{\# F}\cdot
P_F\\
&=\sumsub{\{1\}\cup F\in\calFn\\F\subseteq J\vartriangle (J+1)}2^{1+\# F}\cdot
P_{\{1\}\cup F}\,.
\end{align*}
Thus, to obtain the expressions in~\eqref{E:Dtopeaks-D} for $Y_{\{1\}\cup J}$ and $Y_{\{1'\}\cup J}$,
we only need to show that
\[\psi(Y_{\{1\}\cup J})=\psi(Y_{\{1'\}\cup J})\,.\]
This may be seen as follows. Consider the involution $\rho:\calDn\to\calDn$ given by right
multiplication with the element $\bar{1}\bar{2}3\ldots n$. In other words, $\rho(w)$ is obtained
from $w$ by reversing the sign of both $w_1$ and $w_2$. Note that $w$ has a descent at $1$
but not at $1'$ if and only if $w_1>\abs{w_2}$, while $w$ has a descent at $1'$ but not
at $1$ if and only if $-w_1>\abs{w_2}$. Therefore, $\rho$ restricts to a bijection
between $\{w\in\calDn\, \mid\, \Des(w)=\{1\}\cup J\}$ and $\{w\in\calDn\, \mid\, \Des(w)=\{1'\}\cup J\}$,
for any $J\subseteq\{2,\ldots,n-1\}$. Also, since $\psi$ simply forgets the signs, $\psi\rho=\psi$.
Therefore, $\psi(Y_{\{1\}\cup J})=\psi(Y_{\{1'\}\cup J})$. This completes the proof
of~\eqref{E:Dtopeaks-D}.
Formulas~\eqref{E:Dtopeaks-B} can be obtained similarly, by making use of
Proposition~\ref{P:Btopeaks-B}.
\epf

\section{The peak algebra}\label{S:peakalgebra}

\bd\label{D:peakalgebra} The peak algebra $\pppn$ is the
subspace of the descent algebra $\SolA$ linearly spanned by the elements
$\{P_F\}_{F\in\calFn}$ defined by~\eqref{E:defpeakone} or~\eqref{E:defpeaktwo}.
\ed

We will show below that $\pppn$ is indeed a unital subalgebra of $\SolA$.
As explained in the Introduction, this differs from the object considered
in~\cite{Nym,Sch}.
The connection between the two will be clarified in
Section~\ref{S:ideals}.

Since the elements $P_F$ are sums over disjoint
classes of permutations, they are linearly independent.
Thus, they form a basis of $\pppn$ and
\[\dim\pppn=f_n\,,\]
 the $n$-th Fibonacci number as in~\eqref{E:fibonacci}.

The following is one of the main results of this paper.
\bt\label{T:peakalgebra} $\pppn$ is a
subalgebra of the descent algebra of the symmetric group $\calSn$. Moreover, it
coincides with the image of the map $\varphi: \SolB\to\SolA$ and also
with that of $\psi: \SolD\to\SolA$.
\et
\bpf Since $\varphi$ and $\psi$ are morphisms of algebras, it is enough to
show that their images coincide with $\pppn$. We already know
that
\[\im(\varphi)\subseteq\im(\psi)\subseteq\pppn\]
 from Proposition~\ref{P:Dtopeaks} and Corollary~\ref{C:BDA}. We will show that $\varphi:\SolB\to\pppn$
is in fact surjective.

Define an order on the subsets in $\calFn$
by declaring that  $E<F$ if the maximum element in the
set $E \vartriangle F$ is an element of $F$. According to Lemma~\ref{L:total} (below), this a total order.

For each $E \in \calFn$ consider the shifted set $E-1 \subseteq \Bsets$.
 We first show that $P_F$ does not appear in the sum
$\varphi(X_{E-1})$ for any $E<F$.  Let $a$ be the maximum element
in $E \vartriangle F$.  By the choice of the total order, $a \in F$ and $a \notin E$.  Suppose $a \in
E-1$. Then $a+1 \in E$, and since $a$ is the maximum element of $E
\vartriangle F$, we would have $a+1 \in F$.  But then both $a$ and $a+1\in F$, which contradicts the
fact that $F\in\calFn$.  Thus, $a\notin E-1$, and hence $F \not \subseteq (E-1)\cup E$.  It follows
from~\eqref{E:Btopeaks-B} that
$P_F$ does not appear in the sum $\varphi(X_{E-1})$ for any $E \in \calFn$, $E<F$, while it does appear
in $\varphi(X_{F-1})$. In other words, the matrix of the restriction of $\varphi$ to the subspace of
$\SolB$ spanned by the elements $\{X_{F-1}\,\mid\, F \in \calFn\}$ is unitriangular (with respect
to the chosen total order). In particular, $\varphi$ is surjective.
\epf

\bl\label{L:total}
The relation $E<F\iff \max(E\vartriangle F)\in F$ defines a total order on the finite subsets of
any fixed totally ordered set.
\el
\bpf Let $X$ be the totally ordered set. Encode subsets $E$ of $X$ as sequences
$x_E\in\Z_2^X$ in the standard way. Then $x_{E\vartriangle F}=x_E+x_F$. Therefore,
the relation $E<F$ corresponds to the right lexicographic order on $\Z_2^X$ (defined from
$0<1$ in $\Z_2$), which is clearly total.
\epf

Since the elements $P_F$ consist of sums over disjoint classes of permutations, the structure constants
for the multiplication of $\pppn$ on this basis are necessarily non-negative.
Tables~\ref{T:P2}-\ref{T:P4} describe the multiplication for the first peak algebras. The structure
constants for the product
$P_F\cdot P_G$ are found at the intersection of the row indexed by
$P_F$ and the column indexed by $P_G$, and they are given in the same order as the rows or columns.
For instance,  in $\ppp{4}$,
\[P_{\{2\}}\cdot P_{\{1\}}=2P_\emptyset+2P_{\{1\}}+2P_{\{2\}}+3P_{\{3\}}+3P_{\{1,3\}}\,.\]
Note that the peak algebras are not commutative in general.

\begin{table}[!h]
\[
\begin{array}{|c|c|c|}\hline
 & P_{\emptyset} & P_{\{1\}}  \\ \hline
P_{\emptyset} & (1,0) & (0,1) \\ \hline
P_{\{1\}} & (0,1) & (1,0)  \\ \hline
\end{array}
\]
\caption{Multiplication table for $\ppp{2}$.}\label{T:P2}
\end{table}

\begin{table}[!h]
\[
\begin{array}{|c|c|c|c|}\hline
 & P_{\emptyset} & P_{\{1\}} & P_{\{2\}} \\ \hline
P_{\emptyset} & (1,0,0) & (0,1,0)& (0,0,1) \\ \hline
P_{\{1\}} & (0,1,0) & (2,1,2) & (1,1,1) \\ \hline
P_{\{2\}} & (0,0,1) & (1,1,1) & (1,1,0) \\ \hline
\end{array}
\]
\label{T:P3}
\caption{Multiplication table for $\ppp{3}$.}
\end{table}

\begin{table}[!h]
\[
\begin{array}{|c|c|c|c|c|c|} \hline
 & P_{\emptyset} & P_{\{1\}}&P_{\{2\}} & P_{\{3\}} & P_{\{1,3\}} \\  \hline
P_{\emptyset} & (1,0,0,0,0) & (0,1,0,0,0) & (0,0,1,0,0) &(0,0,0,1,0) &
(0,0,0,0,1) \\ \hline
P_{\{1\}} & (0,1,0,0,0) & (3,2,2,2,2) & (2,2,3,2,2)& (1,1,0,1,2) &
(1,1,2,2,1) \\ \hline
P_{\{2\}} & (0,0,1,0,0) & (2,2,2,3,3) & (3,3,2,3,3)& (1,1,1,1,1) &
(2,2,2,1,1) \\ \hline
P_{\{3\}} & (0,0,0,1,0) & (1,1,1,0,1) & (1,1,1,1,1)& (1,0,1,0,0) &
(0,1,0,1,1) \\ \hline
P_{\{1,3\}} & (0,0,0,0,1) & (1,1,2,2,1) & (2,2,1,2,2) & (0,1,1,0,0) & (2,1,1,1,1)\\ \hline
\end{array}
\]
\label{T:P4}
\caption{Multiplication table for $\ppp{4}$.}
\end{table}

\begin{rem} The proof of Theorem~\ref{T:peakalgebra} shows that the set
$\{\varphi(X_{F-1})\}_{F\in\calFn}$ is a basis of $\pppn$. The structure constants for the product of $\SolB$ on the
basis $\{X_J\}_{J\subseteq\Bsets}$ are non-negative and in fact admit an explicit
combinatorial description~\cite{BB92a}. However, this is no longer true for the basis
$\{\varphi(X_{F-1})\}_{F\in \calFn}$ of $\pppn$. The first counterexample occurs for $n=5$:
\[\varphi(X_{\{2\}})^2=2\varphi(X_{\{2\}})+4\varphi(X_{\{3\}})-2\varphi(X_{\{0,3\}})+14
\varphi(X_{\{1,3\}})\,.\]
We do not know of a basis of $\pppn$ for which the structure constants of the product admit
an explicit combinatorial description.
\end{rem}

\section{Interior peaks and canonical ideals in types B and D}\label{S:ideals}

\subsection{The canonical ideal in type B}\label{S:ideal-B}
\bd\label{D:beta} Let $\beta:\SolB\to\Sol{B_{n-1}}$ be defined by
\begin{equation}\label{E:defbeta-B}
X_J \mapsto \begin{cases} X_{J-1} & \text{ if }0\notin J\\
                          0 & \text{ if }0\in J\,,
             \end{cases}
\end{equation}
for any $J\subseteq\Bsets$. The canonical ideal in type B is $I_n^0=\ker(\beta)$.
\ed
Note that $\beta$ may be equivalently defined by
\begin{equation}\label{E:defbeta-D}
Y_J \mapsto \begin{cases} Y_{J-1} & \text{ if }0\notin J\\
                          -Y_{(J\setminus\{0\})-1} & \text{ if }0\in J\,.
             \end{cases}
\end{equation}
It follows that
\begin{equation}\label{E:defBideals}
I^{0}_n=\Span\{X^0_J\ \mid\  J\subseteq\Asets\}=\Span\{Y^0_J\ \mid\ J\subseteq\Asets\}\,,
\end{equation}
where we have set
\begin{equation}\label{E:defB0D0}
 X^0_J=X_{\{0\}\cup J} \text{ \ and \ } Y^0_J=Y_{\{0\}\cup J}+Y_J\,.
\end{equation}

\bp\label{P:beta} The map $\beta:\SolB\to\Sol{X_{n-1}}$ is a surjective morphism of
algebras. In particular, $I_n^0$ is a two-sided ideal of $\SolB$.
\ep
\bpf The fact that $I_n^0$ is a two-sided ideal of $\SolB$ is proven in~\cite[Theorem 3]{BB92a}
(we warn the reader that the proof of that theorem contains some misprints).
The same argument, which is based on the multiplication rule
for the descent algebra of type B given in~\cite[Theorem 1]{BB92a},
gives the stronger assertion that $\beta$ is a morphism of algebras. Surjectivity is obvious
from~\eqref{E:defbeta-B}.
\epf

\subsection{The ideal of interior peaks}\label{S:intpeaks}
\bd \label{D:intpeaks}
Let $w\in \calSn$. The set of {\em interior peaks} of $w$ is
\[\Peakint(w)=\{i\in\{2,\ldots,n-1\}\ \mid\ w_{i-1}<w_i>w_{i+1}\}\,.\]
\ed

Thus,
\begin{equation}\label{E:intpeaks}
\Peak(w)=\begin{cases}\{1\}\cup\Peakint(w) & \text{ if }1\in\Des(w),\\
                              \Peakint(w) & \text{ if not.}\rule{0pt}{18pt}
                              \end{cases}
\end{equation}
Let $\calFnint$ be the collection of those $F\in\calFn$ for which $1\notin F$.
For any $w\in \calSn$, $\Peakint(w)\in\calFnint$, and it is easy to see that any
$F\in\calFnint$ can be realized in this way.

For each $F\in\calFnint$, define an element of the group algebra of $\calSn$ by
\begin{equation}\label{E:defpeakintone}
\Pint_F:=\sum_{\Peakint(w)=F}w\,.
\end{equation}
These are the elements considered in Nyman's original work~\cite{Nym} and also in
Schocker's~\cite{Sch}. The operator
$\Lambdaint(J):=\{i\in J\ \mid\ i\neq 1,\ i-1\notin J\}$ satisfies
$\Lambdaint(\Des(w))=\Peakint(w)$; thus, we also have
\begin{equation}\label{E:defpeakinttwo}
\Pint_F:=\sum_{\Lambdaint(J)=F}Y_J\,.
\end{equation}

Note that if $F\in\calFnint$ and $2\notin F$ then $\{1\}\cup F\in\calFn$.
\bl\label{L:peakvspeakint}For any $F\in\calFnint$,
\begin{equation}\label{E:peakvspeakint}
\Pint_F=\begin{cases} P_F & \text{ if }2\in F\\
P_F+P_{\{1\}\cup F} & \text{ if }2\notin F\,.\end{cases}
\end{equation}
\el
\bpf Consider first the case when $2\in F$. Then, for any $w\in \calSn$ with
$\Peakint(w)=F$, $w_1<w_2>w_3$; in particular, $1\notin\Des(w)$. Hence
$\Peak(w)=\Peakint(w)$ and it follows from~\eqref{E:defpeakone}
and~\eqref{E:defpeakintone} that $P_F=\Pint_F$.

If $2\notin F$, then the class of permutations $w\in \calSn$ with $\Peakint(w)=F$
splits into two classes, according to whether they have a descent at $1$ or not.
In view of~\eqref{E:intpeaks}, the first class consists of all those $w$ for
which $\Peak(w)=\{1\}\cup F$, and the second of those for which $\Peak(w)=F$;
whence $\Pint_F=P_{\{1\}\cup F}+P_F$.
\epf

\bd\label{D:peakideal} The peak ideal $\pppnint$ is the subspace of $\pppn$ linearly
spanned by the elements $\{\Pint_F\}_{F\in\calFnint}$.
\ed

We will show below that $\pppnint$ is indeed a two-sided ideal of $\pppn$. In particular, it
follows that $\pppnint$ is closed under the product of $\SolA$. In ~\cite{Nym,Sch},
$\pppnint$ is regarded as a non
unital subalgebra of $\SolA$ (and it is referred to as the peak algebra). The present work reveals the
existence of the larger unital subalgebra $\pppn$ of $\SolA$, as well as the extra structure of
 $\pppnint$ as an ideal of $\pppn$.

 First, consider the map $\pi:\pppn\to\ppp{n-2}$ defined by
\begin{equation}\label{E:defpi}
P_F\mapsto
\begin{cases}
P_{F-2} & \text{ if neither $1$ nor $2$ belong to $F$}\\
-P_{F\setminus\{1\}-2} & \text{ if }1\in F\\
0 & \text{ if }2\in F
\end{cases}
\end{equation}
for any $F\in\calFn$. Note that~\eqref{E:defpi} covers all possibilities, since
$1$ and $2$ cannot belong to $F$ simultaneously. For the same reason, if $1\in F$,
the set $F\setminus\{1\}-2$ belongs to $\calF_{n-2}$.

As we will see below (Remark \ref{R:nomorphism}), the map $\beta:\SolB\to\Sol{B_{n-1}}$
 does not descend to peak algebras. However, the composite
\begin{equation}\label{E:defbeta-B2}
\beta^2:\SolB\map{\beta}\Sol{B_{n-1}}\map{\beta}\Sol{B_{n-2}}\,, \qquad
X_J \mapsto \begin{cases} X_{J-2} & \text{ if }0,1\notin J\\
                          0 & \text{ if $0$ or $1\in J$}
             \end{cases}
\end{equation}
does.  Indeed, it descends precisely to the map $\pi$
of~\eqref{E:defpi}.

\bp\label{P:pi} The map $\pi:\pppn\to\ppp{n-2}$ is a surjective morphism of
algebras. Moreover, there is a commutative diagram
\[\xymatrix{ {\SolB}\ar[r]^{\beta^2}\ar[d]_{\varphi} &
{\Sol{B_{n-2}}}\ar[d]^{\varphi}\\
{\pppn}\ar[r]_{\pi}&{\ppp{n-2}} }\]
\ep
\bpf  It suffices to show that the diagram commutes. That $\pi$ is a surjective
morphism of algebras follows, since so are $\varphi$ and $\beta$ (Theorem~\ref{T:peakalgebra}
and Proposition~\ref{P:beta}).

Take $J\subseteq\Bsets$. Suppose first that $0$ or $1\in J$, so that $\beta^2(X_J)=0$. In this case,
$1\in J\cup(J+1)$ and hence
\begin{align*}
\pi\varphi(X_J) &\ \equal{\eqref{E:Btopeaks-B}}\
 2^{\#J}\cdot\!\!\!\!\sumsub{F\in\calFn\\ F\subseteq J\cup(J+1)} \pi(P_F)\\
&\ \equal{\eqref{E:defpi}}\ 2^{\#J}\cdot\!\!\!\!\sumsub{F\in\calFn\\ F\subseteq J\cup(J+1)\\1\notin F,
2\notin F}P_{F-2}-2^{\#J}\cdot\!\!\!\!\sumsub{F\in\calFn\\ F\subseteq J\cup(J+1)\\1\in F,
2\notin F}P_{F\setminus\{1\}-2}\\
&\ =\ 2^{\#J}\cdot\!\!\!\!\sumsub{F\in\calFn\\ F\subseteq J\cup(J+1)\\1\notin F,
2\notin F}P_{F-2}-2^{\#J}\cdot\!\!\!\!\sumsub{F'\in\calFn\\ F'\subseteq J\cup(J+1)\\1\notin F',
2\notin F'}P_{F'-2}=0\,.  \rule{0pt}{18pt}
\end{align*}
Consider now the case when $0,1\notin J$. In this case,
\[\varphi\beta^2(X_J)\ =\ \varphi(X_{J-2})\ \equal{\eqref{E:Btopeaks-B}}\ 2^{\#J}\cdot\!\!\!\!\sumsub{F\in\calFn\\
F\subseteq (J-2)\cup(J-1)} P_F\,,\]
since the cardinal of the shifted set $J-2$ coincides with the cardinal of $J$.
On the other hand, since $1\notin J\cup(J+1)$, we have
\[\pi\varphi(X_J)\ \equal{\eqref{E:Btopeaks-B}}\ 2^{\#J}\cdot\!\!\!\!\sumsub{F\in\calFn\\
F\subseteq J\cup(J+1)} \pi(P_F)
\ \equal{\eqref{E:defpi}}\ 2^{\#J}\cdot\!\!\!\!\sumsub{F\in\calFn\\
F\subseteq J\cup(J+1)\\1,2\notin F}\!\! P_{F-2}\ =\ 2^{\#J}\cdot\!\!\!\!\sumsub{F\in\calFn\\
F\subseteq (J-2)\cup(J-1)}\!\!\! P_{F}\,.\]
Thus, $\varphi\beta^2(X_J)=\pi\varphi(X_J)$ and the proof is complete.
\epf

Since the elements $\Pint_F$ are sums over disjoint classes of permutations,
they are linearly independent. Thus, they form a basis of $\pppnint$ and
\[\dim\pppnint=\#\calFnint=\#\calF_{n-1}=f_{n-1}\,,\]
 the $(n-1)$-th Fibonacci number as in~\eqref{E:fibonacci}.

We can now derive the main result of the section. It may be viewed as an algebraic counterpart
of the recursion $f_n=f_{n-1}+f_{n-2}$.
\bt\label{T:peakideal} $\pppnint$ is a two-sided ideal of $\pppn$. Moreover, for $n\geq 2$,
\[\pppn/\pppnint\cong \ppp{n-2}\,.\]
\et
(For $n\leq1$, we have $\ppp{0}=\ppp{1}=\pppint{1}=\Q$.)
\bpf We know from Proposition~\ref{P:pi} that $\pi$ is  surjective. Therefore
\[\pppn/\ker(\pi)\cong \ppp{n-2}\,.\]
To complete the proof of the theorem it suffices to show that
\[\ker(\pi)=\pppnint\,.\]
Since $\dim\pppn=f_n$, we have
\[\dim\ker(\pi)=f_n-f_{n-2}=f_{n-1}=\dim\pppnint\,.\]
On the other hand, from~\eqref{E:peakvspeakint} and~\eqref{E:defpi} we immediately
see that $\pppnint\subseteq\ker(\pi)$.
\epf

\subsection{From the canonical ideal to the peak ideal}\label{S:Btointpeaks}

 Let $I_n^{0,1}$ denote the kernel of the map
 $\beta^2$. Thus, $I_n^{0,1}$ is another ideal
 of $\SolB$ and $ I_n^0\subseteq I_n^{0,1}$.
 Moreover, from~\eqref{E:defbeta-B} and~\eqref{E:defbeta-B2} we see that
\begin{equation}
I^{0,1}_n=\Span\{X_J\ \mid\ 0 \text{ or }1\in J,\ J\subseteq\Bsets\}\,.
\end{equation}
We will show below that the ideal $I_n^0$ and the larger ideal $I_n^{0,1}$ both map onto $\pppnint$.
First, we provide an explicit formula for the map $\varphi$ in terms of the bases
$\{X^0_J\}_{J\subseteq\Asets}$ and  $\{Y^0_J\}_{J\subseteq\Asets}$
of $I_n^0$ and $\{\Pint_F\}_{F\in\calFnint}$ of $\pppnint$. These turn
out to be completely analogous to the formulas in Propositions~\ref{P:Btopeaks-D}
and~\ref{P:Btopeaks-B}.

\bp\label{P:Btointpeaks} For any $J\in\Asets$,
\begin{equation}\label{E:Btointpeaks}
\varphi(X^0_J)=2^{1+\#J}\cdot\!\!\!\!\sumsub{F\in\calFnint\\F\subseteq J\cup(J+1)}\Pint_{F}
\text{ \ \ and \ \ }
\varphi(Y^0_J)=\!\!\!\!\sumsub{F\in\calFnint\\F\subseteq J\vartriangle(J+1)}\!\!\!
2^{1+\#F}\cdot\Pint_{F}\,.
\end{equation}
\ep
\bpf By
Proposition~\ref{P:Btopeaks-B} we have
\begin{align*}
\frac{1}{2^{1+\#J}}\cdot\varphi(X_{\{0\}\cup J}) &=
\sumsub{F\in\calFn\\F\subseteq\{1\}\cup J\cup(J+1)} P_F\\
&=\sumsub{1\notin F\in\calFn\\2\in F\subseteq J\cup(J+1)} P_F
+\sumsub{1\in F\in\calFn\\2\notin F\subseteq\{1\}\cup J\cup(J+1)} P_F
+\sumsub{1\notin F\in\calFn\\2\notin F\subseteq\{1\}\cup J\cup(J+1)} P_F\\
&\equal{\eqref{E:peakvspeakint}}\sumsub{F\in\calFnint\\2\in F\subseteq J\cup(J+1)} \Pint_F
+\sumsub{F'\in\calFn\\1,2\notin F'\subseteq J\cup(J+1)} P_{\{1\}\cup F'}
+\sumsub{ F\in\calFn\\1,2\notin F\subseteq J\cup(J+1)} P_F\\
&=\sumsub{F\in\calFnint\\2\in F\subseteq J\cup(J+1)} \Pint_F
+\sumsub{F\in\calFn\\1,2\notin F\subseteq J\cup(J+1)} P_{\{1\}\cup F}
+ P_F\\
&\equal{\eqref{E:peakvspeakint}}\sumsub{F\in\calFnint\\2\in F\subseteq J\cup(J+1)} \Pint_F
+\sumsub{F\in\calFnint\\2\notin F\subseteq J\cup(J+1)}\Pint_{F}\\
&=\sumsub{F\in\calFnint\\F\subseteq J\cup(J+1)}\Pint_{F}\,.
\end{align*}
The formula for $Y^0_J$ can be similarly deduced from Proposition~\ref{P:Btopeaks-D}, or
from the formula for $X^0_J$, as in the proof of Proposition~\ref{P:Btopeaks-B}.
\epf

\bt\label{T:Bideals} $\varphi(I_n^{0,1})=\varphi(I_n^0)=\pppnint$.
\et
\bpf  We already know that
$\varphi(I_n^0)\subseteq \varphi(I_n^{0,1})\subseteq\pppnint$. In fact,
the first inclusion is trivial and the second follows from Proposition~\ref{P:pi}, since
$I_n^{0,1}=\ker(\beta^2)$ and $\pppnint=\ker(\pi)$, as seen in the proof of Theorem~\ref{T:peakideal}.

To complete the proof, we show that $\varphi:I_n^0\to\pppnint$ is surjective.
For each $E \in \calFnint$, consider the element  $X^0_{E-1} \in I_n^{0}$.
By Proposition~\ref{P:Btointpeaks},
\[\varphi(X^0_{E-1}) =2^{1+\#E}\cdot\sumsub{F\in\calFnint\\F\subseteq (E-1)
\cup E}\Pint_F\,.\]
As in the proof of Theorem~\ref{T:peakalgebra}, it follows that  $\Pint_F$ does not appear in
$\varphi(X^0_{E-1})$ for any $E < F$, where $E<F$ denotes the total order of
Lemma~\ref{L:total}.  On the other hand,  $\Pint_F$ does appear in $\varphi(X^0_{F-1})$. Therefore, the matrix of the map $\varphi$ restricted to the subspace of $I_n^0$ spanned
by the elements  $X^0_{E-1}$ is unitriangular. Hence, $\varphi:I_n^0\to\pppnint$ is
surjective.
\epf

The following corollary summarizes the results that relate the peak algebra
and the peak ideal to the descent algebras of type B.
\bc\label{C:Bexactsequence} There is a commutative diagram
\[\xymatrix{0\ar[r]& {I_n^{0,1}}\xyinc[r]\ar[d]_{\varphi}  &
{\SolB}\ar[r]^{\beta^2}\ar[d]^{\varphi} & {\Sol{B_{n-2}}}\ar[d]^{\varphi}\ar[r] & 0\\
0\ar[r]& {\pppnint}\xyinc[r] & {\pppn}\ar[r]_{\pi}&{\ppp{n-2}}\ar[r] & 0 }\]
where the rows are exact and the vertical maps are surjective.
\ec
\bpf The surjectivity of the vertical maps was obtained in Theorems~\ref{T:peakalgebra} and~\ref{T:Bideals}.
The exactness of the first row is clear from the definitions of $\beta$ and
$I_n^{0,1}$.
In the course of the proof of Theorem~\ref{T:peakideal} we obtained that
$\ker(\pi)=\pppnint$.  Finally,
Proposition~\ref{P:pi} gave the surjectivity of $\pi$ and  the commutativity
of the  square.
\epf

\begin{rem}\label{R:nomorphism}
There is no morphism of algebras $\pppn\to\ppp{n-1}$ fitting in
a commutative diagram
\[\xymatrix{ {\SolB}\ar[r]^{\beta}\ar[d]_{\varphi} & {\Sol{B_{n-1}}}\ar[d]^{\varphi}\\
{\pppn}\ar@{-->}[r]&{\ppp{n-1}} }\]

For, such a map would be surjective, since so are $\varphi$ and $\beta$, and hence
its kernel would have dimension $\dim\pppn-\dim\ppp{n-1}=f_{n-2}$. But this kernel
must contain
\[\varphi(\ker(\beta))=\varphi(I_n^0)=\pppnint \text{ (by Theorem~\ref{T:Bideals})}\]
which has dimension $f_{n-1}$, a contradiction.
\end{rem}

\subsection{The canonical ideal in type D}\label{S:ideal-D}

There are similar results relating peaks to the descent algebras of type D.
We present them in this section.

Consider the map $\gamma:\SolD\to\Sol{B_{n-2}}$ defined by
\begin{equation}\label{E:defgamma}
X_J\mapsto\begin{cases}X_{J-2} & \text{ if }1,1'\notin J\\
0 & \text{ if $1$ or $1'\in J$.}  \end{cases}
\end{equation}
It follows from the proof of Theorem 2 in~\cite{BV} that $\gamma$ is a morphism
of algebras.  Moreover, it is clear from equations~\eqref{E:BDmap2}, \eqref{E:defbeta-B} and \eqref{E:defgamma}
that the diagram
\begin{equation}\label{E:gammabeta}
\xymatrix@R=1.2pc@C=0.6pc{{\SolB}\ar[rr]^{\chi}\ar[dr]_{\beta^2} & &{\SolD}\ar[ld]^{\gamma}\\
& {\Sol{B_{n-2}}} }
\end{equation}
commutes.

Let  $I_n^{1',1}$ denote the kernel of the map $\gamma$, an ideal
 of $\SolD$. From~\eqref{E:defgamma}  we see that
\begin{equation}\label{E:defDideal}
I^{1',1}_n=\Span\{X_J\ \mid\ 1' \text{ or }1\in J,\ J\in\Dsets\}\,.
\end{equation}

The following theorem describes the relationship between the peak algebra
and the peak ideal to the descent algebras of type D and B.
\bt\label{T:Dexactsequence} There is a commutative diagram
\[\xymatrix{0\ar[r]& {I_n^{1',1}}\xyinc[r]\ar[d]_{\psi}  &
{\SolD}\ar[r]^{\gamma}\ar[d]^{\psi} & {\Sol{B_{n-2}}}\ar[d]^{\varphi}\ar[r] & 0\\
0\ar[r]& {\pppnint}\xyinc[r] & {\pppn}\ar[r]_{\pi}&{\ppp{n-2}}\ar[r] & 0 }\]
where the rows are exact and the vertical maps are surjective.
\et
\bpf The exactness of the first row is clear from the definitions of $\gamma$ and
$I_n^{1',1}$. The second row is exact according to Corollary~\ref{C:Bexactsequence}.

The commutativity of the square
\[\xymatrix{ {\SolD}\ar[r]^{\gamma}\ar[d]_{\psi} & {\Sol{B_{n-2}}}\ar[d]^{\varphi}\\
{\pppn}\ar[r]_{\pi}&{\ppp{n-2}} }\]
can be shown in the same way as in Proposition~\ref{P:pi}. Together
with~\eqref{E:gammabeta}, this allows us to construct a commutative diagram
 \[\xymatrix{0\ar[r]& {I_n^{0,1}}\xyinc[r]\ar[d]_{\chi}  &
{\SolB}\ar[rd]^{\beta^2}\ar[d]^{\chi} & &\\
0\ar[r]& {I_n^{1',1}}\xyinc[r]\ar[d]_{\psi}  &
{\SolD}\ar[r]^{\gamma}\ar[d]^{\psi} & {\Sol{B_{n-2}}}\ar[d]^{\varphi}\ar[r] & 0\\
0\ar[r]& {\pppnint}\xyinc[r] & {\pppn}\ar[r]_{\pi}&{\ppp{n-2}}\ar[r] & 0 }\]
Since by Theorem~\ref{T:Bideals} the map $I_n^{0,1}\to\pppnint$ is surjective, so must be the
map $I_n^{1',1}\to\pppnint$. This completes the proof.
\epf

\section{Commutative subalgebras of the descent and peak algebras} \label{S:commutative}

Grouping permutations according to the number of descents leads to the following elements of the descent algebra of type A:
\[ y_j := \sumsub{u\in S_n \\\#\Des(u)=j} u=\sumsub{J\subseteq\Asets\\\#J=j}Y_J\,,\]
for $j=0,\ldots, n-1$. It was shown by Loday that the linear span of these elements forms
a commutative subalgebra of $\SolA$~\cite[Section 1]{Lod89}, see also \cite[Remark 4.2]{GReu}.
Moreover, this subalgebra is generated by $y_1$ and is semisimple.
Let us denote it by $\sA$.

In this section we derive an analogous result for the linear span of the
sums of permutations with a given number of peaks. Once again, it is possible to easily derive
this fact from the analogous result for type B.

\subsection{Commutative subalgebras of the descent algebra of type B}\label{S:comm-B}
Consider the following elements of the descent algebra of type B:
\begin{equation}\label{E:solB-d}
y_j :=\sumsub{w \in B_n \\ \#\Des(w)=j}\!\!\!w\ =\sumsub{J\subseteq\Bsets\\\#J=j}Y_J\,,
\end{equation}
for each $j=0,\ldots,n$. Let $\sB$ denote the subspace of $\SolB$ linearly spanned by these $n{+}1$ elements.

Let also
\begin{equation}\label{E:solB-b}
x_j :=\sumsub{J\subseteq\Bsets\\\#J=j}X_J\,,
\end{equation}
for $j=0,\ldots,n$. These elements also span $\sB$, since it is easy to see that
\begin{equation}\label{E:btod}
x_j=\sum_{i=0}^j\binom{n-i}{j-i}\cdot y_i\,.
\end{equation}
Bergeron and Bergeron showed that $\sB$ is
a commutative semisimple subalgebra of $\SolB$~\cite[Section 4]{BB92b}.
Mahajan showed that one can actually do better~\cite{Mah01}. Below, we will formulate his results
in our notation (Mahajan's work is in the context of random walks on chambers of
Coxeter complexes), together with a small extension.
To this end,  recall the bases $\{Y^0_J\}_{J\subseteq\Asets}$ and
$\{X^0_J\}_{J\subseteq\Asets}$ of the ideal $I_n^0$ of $\SolB$ from~\eqref{E:defB0D0},
and consider the elements
\begin{equation}\label{E:in0}
y^0_j:=\sumsub{J\subseteq\Asets\\\#J=j-1}Y^0_J \text{ \ \ and \ \ }
x^0_j:=\sumsub{J\subseteq\Asets\\\#J=j-1}X^0_J=\sum_{i=1}^{j}\binom{n-i}{j-i}\cdot y^0_i
\end{equation}
for $j=1,\ldots,n$. Let $i_n^0$ denote the subspace of $I_n^0$ linearly
spanned by (either of) these two sets of $n$ elements and set
\[\hsB:=\sB+i_n^0\,,\]
the subspace of $\SolB$ linearly spanned by $y_0,\ldots,y_n,y_1^0,\ldots,y_n^0$.
It is easy to see that the only relation among these elements is
\begin{equation}\label{E:Bsumall}
\sum_{j=0}^n y_j=\sum_{w\in\calBn}w=\sum_{j=1}^{n} y_j^0\,, \text{ \ or in other words, \ }x_n=x_n^0\,;
\end{equation}
thus, $\sB\cap i_n^0=\Span\{x_n\}$ and $\dim\hsB=2n$.

\bt\label{T:solB}
$\hsB$ is a commutative semisimple subalgebra of $\SolB$. Moreover:
\begin{itemize}
\item[(a)] $\hsB$ is generated as an algebra by $y_1$ and $y_1^0$,
\item[(b)] $\sB$ is the subalgebra of $\hsB$ generated by $y_1$, it is also semisimple,
\item[(c)] $i_n^0$ is the two-sided ideal of $\hsB$ generated by $y_1^0$,
\item[(d)] $\dim\hsB=2n$,  $\dim\sB=n+1$,  $\dim i_n^0=n$.
\end{itemize}
\et
\bpf  All these assertions are due to Mahajan~\cite[Section 4.2]{Mah01}, except for (c).
Since $i_n^0\subseteq I_n^0$ and $\hsB=\sB+i_n^0$, we have
\[I_n^0\cap\hsB=(I_n^0\cap\sB)+i_n^0\,.\]
By Definition~\ref{D:beta}, $I_n^0=\ker(\beta)$. Hence,
 $I_n^0\cap\sB=\ker(\beta|_{\sB})$. Now, from our analysis of $\beta$
 in Proposition~\ref{P:solbeta} (below), this kernel is $1$-dimensional, spanned
 by the element $\sum_{i=0}^n y_i$, which belongs to $i_n^0$ by~\eqref{E:Bsumall}.
 Therefore,
 \[I_n^0\cap\hsB=i_n^0\,.\]
Since $I_n^0$ is an ideal of $\SolB$, it follows that $i_n^0$ is an ideal of $\hsB$.
Mahajan also shows that $i_n^0$ is generated by $y_1^0$ as a non-unital algebra, hence also as
an ideal.
\epf

\subsection{From number of descents to number of peaks}\label{S:sBtowp}

The maximum number of peaks of a permutation $u\in\calSn$ is $\ipartn$.
For each $j=0,\ldots,\ipartn$, define an element of the peak algebra $\pppn$ by
\begin{equation}\label{E:wp}
p_j\ := \!\!\sumsub{w \in \calSn\\ \#\Peak(w)=j}\!\!\! w\ =\ \sumsub{F\in\calFn\\\#F=j}P_F\,.
\end{equation}
It is somewhat surprising that this naive analog of~\eqref{E:solB-d} is indeed  the
right one, given that the expression for $\varphi(Y_J)$ in~\eqref{E:Btopeaks-D} involves subsets $F$
of various cardinalities. However, we have:

\bp\label{P:sBtowp}
The restriction of $\varphi:\SolB\to\pppn$ to $\sB$  is explicitly given by
\begin{equation} \label{E:sBtowp-d}
y_j  \mapsto \sum_{i=0}^{\min(j,n-j)} 2^{2i} \binom{n-2i}{j-i}\cdot p_i
\end{equation}
for $j=0,\ldots,n$.
\ep
\bpf Fix $j\in\{0,\ldots,n\}$. By~\eqref{E:solB-d} and~\eqref{E:Btopeaks-D},
\[\varphi(y_j)=\sumsub{J\subseteq\Bsets\\\#J=j}
\sumsub{F\in\calFn\\F\subseteq J\vartriangle (J+1)}2^{\# F}\cdot P_F\,.\]
To derive~\eqref{E:sBtowp-d} we have to show that, for each $i=0,\ldots,\ipartn$ and $F\in\calFn$
with $\#F=i$, we have
\[\#\{J\subseteq\Bsets\ \mid\ F\subseteq J\vartriangle (J+1)
\text{ and }\#J=j\}=\begin{cases} 2^i\binom{n-2i}{j-i} & \text{ if }i\leq j,\,n-j,\\
0 & \text{ if not.}\end{cases}\]

Let $F=\{s_1, s_2, \ldots, s_i\}$ and consider the
sets $\{s_h-1, s_h\}$.  These sets are disjoint since $F\in\calFn$.
Now, $F\subseteq J\vartriangle(J+1)$ if
and only if $J$ contains exactly one of the elements $s_h-1$, $s_h$
for each $h=1,\ldots,i$. In particular, there is no such set $J$ unless $i\leq j$. In this
case, $J$ is determined by the choice of one element from each set $\{s_h-1, s_h\}$
 plus the choice  of $j-i$ elements from the
remaining $n-2i$ elements in $\Bsets$ not occurring in any of those sets (which is possible
only if $i\leq n-j$).
The total number of choices is thus $2^i\binom{n-2i}{j-i}$.
\epf

\begin{rem} It follows immediately from~\eqref{E:sBtowp-d} that $\varphi(y_j)=\varphi(y_{n-j})$
for each $j=0,\ldots,n$. This is also a consequence of Remark~\ref{R:kerphi}.
\end{rem}

The maximum number of interior peaks of a permutation $u\in\calSn$ is $\ipart{n-1}$.
For each $j=1,\ldots,\ipart{n+1}$
define an element of the peak ideal by
\begin{equation}\label{E:pint}
\pint_j\ :=\!\! \sumsub{w \in \calSn\\ \#\Peakint(w)=j-1}\!\!\!\! w\ =\sumsub{F\in \calFnint\\ \#F=j-1} \Pint_F.
\end{equation}

\bp\label{P:i0topint}
The restriction of $\varphi$ to $i_n^0$ is given by
\begin{equation} \label{E:i0topint-d}
y_j^0  \mapsto \sum_{i=1}^{\min(j,n+1-j)} 2^{2i-1}\binom{n-2i+1}{j-i}\cdot \pint_i
\end{equation}
for each $j=1,\ldots,n$.
\ep
\bpf The same argument as in the proof of Proposition~\ref{P:sBtowp} allows
us to deduce this result from~\eqref{E:Btointpeaks}.
\epf

\begin{rem} We may rewrite the definition of the elements $y^{0}_{j}$ and $\pint_j$  as
\[y^0_j\ =\!\!\sumsub{J\subseteq\Bsets\\\#(J\setminus \{0\})=j-1}\!\!\!Y_J
\text{ \ and \ }\
\pint_j\ =\!\!\sumsub{F\in\calFn\\\#(F\setminus \{1\})=j-1}\!\!\!P_F\,.\]
The former follows from~\eqref{E:in0} and~\eqref{E:defB0D0}, the latter from~\eqref{E:pint} and~\eqref{E:peakvspeakint}.
\end{rem}

\subsection{Commutative subalgebras of the peak algebra}\label{S:comm-peaks}

\bd\label{D:wp}  Let $\wp_n$ denote the subspace of the peak algebra $\pppn$
linearly spanned by the elements $\{p_j\}_{j=0,\ldots,\ipartn}$, let $\wpint_n$ denote the
subspace of the peak ideal $\pppnint$ linearly spanned by the elements
$\{\pint_j\}_{j=1,\ldots,\ipart{n+1}}$, and let
\[\hwp_n:=\wp_n + \wpint_n\,.\]
\ed

\bp\label{P:phionto} The map $\varphi:\SolB\onto\pppn$ maps $\hsB$ onto $\hwp_n$,
$\sB$ onto $\wp_n$ and $i^0_n$ onto $\wpint_n$.
\ep
\bpf  It is enough to prove the last two assertions.
By Proposition~\ref{P:sBtowp},  $\varphi$ maps $\sB$ into $\wp_n$.
Moreover,~\eqref{E:sBtowp-d} also shows that, for
any $j=0,\ldots,\ipartn$,  $p_j$ appears in $\varphi(y_j)$, but not
in $\varphi(y_i)$ for $i<j$. Thus, the matrix of the restriction of $\varphi$
to the subspace of $\sB$ spanned by these elements $p_j$ is unitriangular. Hence,
$\varphi$ maps $\sB$ onto $\wp_n$. That it maps $i^0_n$ onto $\wpint_n$ follows similarly
from ~\eqref{E:i0topint-d}.
\epf

We can now derive our  main result on commutative subalgebras of the peak algebra.
\bt\label{T:wp}$\hwp_n$ is a commutative semisimple subalgebra of $\pppn$. Moreover:
\begin{itemize}
\item[(a)] $\hwp_n$ is generated as an algebra by $p_1$ and $\pint_1$,
\item[(b)] $\wp_n$ is the subalgebra of $\hwp_n$ generated by $p_1$, it is also semisimple,
\item[(c)] $\wpint_n$ is the two-sided ideal of $\hwp_n$ generated by $\pint_1$,
\item[(d)] $\dim\hwp_n=n$, $\dim\wp_n=\ipartn+1$, $\dim\wpint_n=\ipart{n+1}$.
\end{itemize}
\et
\bpf All these statements follow from the corresponding ones in Theorem~\ref{T:solB},
thanks to Proposition~\ref{P:phionto}, except for (d).
Now, it is clear that both sets$\{p_j\}_{j=0,\ldots,\ipartn}$ and
$\{\pint_j\}_{j=1,\ldots,\ipart{n+1}}$
are linearly independent, and that the only relation among these elements is
\[\sum_{j=0}^{\ipartn} p_j=\sum_{u\in\calSn}u=\sum_{j=1}^{\ipart{n+1}} \pint_j\,.\]
(This also follows from Proposition~\ref{P:wppi} below).
The assertions follow, since $n=\ipartn+\ipart{n+1}$.
\epf

The first multiplication tables for the spanning set $\{p_0,\ldots,
p_{\ipartn},\pint_1,\ldots,\pint_{\ipart{n+1}}\}$ of the algebra $\hwp_n$ are provided below.
The first block gives the multiplication within the subalgebra $\wp_n$
and the others describe the multiplication with elements of the ideal $\wpint_n$.
\begin{table}[!h]
\[\begin{array}{|c|c|c||c|}\hline
 & p_0 & p_1 & \pint_1 \\ \hline
p_0 & (1,0,0) & (0,1,0) & (0,0,1)\\ \hline
p_1 & (0,1,0) & (1,0,0) & (0,0,1) \\ \hline \hline
\pint_1 & (0,0,1) & (0,0,1) & (0,0,2)   \\ \hline
\end{array}\]
\caption{Multiplication table for $\hwp_2$.}\label{T:wp2}
\end{table}
\begin{table}[!h]
\[\begin{array}{|c|c|c||c|c|}\hline
 & p_0 & p_1 & \pint_1 &\pint_2\\ \hline
p_0 & (1,0,0,0) & (0,1,0,0) & (0,0,1,0) &(0,0,0,1)\\ \hline
p_1 & (0,1,0,0) & (5,4,0,0) & (0,0,3,4) & (0,0,2,1)\\ \hline \hline
\pint_1 & (0,0,1,0) & (0,0,3,4) & (0,0,3,2) & (0,0,1,2)\\  \hline
\pint_2 & (0,0,0,1) & (0,0,2,1) & (0,0,1,2) & (0,0,1,0)\\  \hline
\end{array}
\]
\caption{Multiplication table for $\hwp_3$.}\label{T:wp3}
\end{table}
\begin{table}[!h]
\[
\begin{array}{|c|c|c|c||c|c|} \hline
 & p_0 & p_1 & p_2  & \pint_1 & \pint_2\\  \hline
p_0 & (1,0,0,0,0) & (0,1,0,0,0) & (0,0,1,0,0) &(0,0,0,1,0)   &(0,0,0,0,1)  \\ \hline
p_1 & (0,1,0,0,0) & (15,13,15,0,0) & (3,4,3,0,0) & (0,0,0,6,6) & (0,0,0,12,12) \\ \hline
p_2 & (0,0,1,0,0) & (3,4,3,0,0) & (2,1,1,0,0) & (0,0,0,1,2)& (0,0,0,4,3)\\ \hline \hline
\pint_1 & (0,0,0,1,0) &(0,0,0,6,6) & (0,0,0,1,2) & (0,0,0,4,2) & (0,0,0,4,6)\\ \hline
\pint_2 & (0,0,0,0,1) &(0,0,0,12,12) & (0,0,0,4,3) & (0,0,0,4,6) & (0,0,0,12,10)\\ \hline
\end{array}
\]
\caption{Multiplication table for $\hwp_4$.}\label{T:wp4}
\end{table}

\br \label{S:otherwork}
\be
\item The elements $\pint_j$ are also considered in recent work of Schocker:
in~\cite[Theorem 9]{Sch}, it is shown that they span a commutative non-unital subalgebra.
Our results (Theorem~\ref{T:wp}) provide a more complete picture: these elements actually span a
two-sided ideal of a
larger commutative unital subalgebra, also defined by grouping permutations according to the
number of (not necessarily interior) peaks.

\item  Doyle and Rockmore have considered two objects closely related to $\wp_n$ and
$\wpint_n$ and related them to
descents of signed permutations~\cite[Sections 5.3-5.5]{DR}. The difference stems
from the fact that in their work permutations are grouped according to the number of
peaks {\em and} valleys
(since they do not consider $n$ as a possible valley, this number is not just
twice the number of peaks).

\item The work of Mahajan~\cite{Mah01} presents a powerful axiomatic construction of
commutative semisimple subalgebras of descent algebras of Coxeter groups.
In particular, several such subalgebras of $\SolB$
are described.

\item More information on the commutative subalgebra $\sB$ of Section~\ref{S:comm-B}
can be found in Chow's thesis~\cite[Sections 4.2, 4.3 and 5.1]{Chow}.

\item Neither the commutative subalgebra $\wp_n$ of $\pppn$ nor the ideal
$\wpint_n$ are contained in the
commutative subalgebra $\sA$ of $\SolA$ discussed in the introduction to Section~\ref{S:commutative}
(counterexamples already exist for low values of $n$). There is, however, a ``correct'' type B analog
of $\sA$: it is a commutative subalgebra of the algebra of Mantaci and Reutenauer
(Section~\ref{S:ManReu}). We will not pursue this point in this paper.
\ee
\er

\subsection{Exact sequences for the commutative subalgebras} \label{S:exact}
Next we consider the morphisms of algebras $\beta$ and $\pi$ (Section~\ref{S:intpeaks}) restricted to
the subalgebras $\sB$ and $\wp_n$. Note that $i_n^0\subseteq I_n^0=\ker(\beta)$ and
$\wpint_n\subseteq\pppnint=\ker(\pi)$.

\bp\label{P:solbeta} For any $j=0,\ldots,n$,
\begin{equation}\label{E:solbeta}
\beta(y_j)= \begin{cases} y_0 & \text{ if }j=0\,,\\
                          y_j-y_{j-1} & \text{ if }0<j<n\,,\\
                          -y_{n-1} & \text{ if } j=n\,;
             \end{cases}  \text{ \ \ or equivalently, \ \ }
\beta(x_j)= \begin{cases} x_j & \text{ if } 0\leq j<n\,,\\
                          0 & \text{ if } j=n\,.
             \end{cases}\end{equation}
Therefore, $\beta$ maps $\sB$ onto $\sol{B_{n-1}}$ with a $1$-dimensional kernel spanned by
\[x_n=\sum_{j=0}^n y_j=\sum_{w\in\calBn}w\,.\]
\ep
\bpf Formulas~\eqref{E:solbeta} follow easily from~\eqref{E:defbeta-B} and ~\eqref{E:defbeta-D};
the rest of the assertions follow at once.
\epf

Even though $\varphi:\sB\to\wp_n$ does not map the $y_j$ basis onto the $p_j$ basis,
the result for the restriction to $\pi$ is completely
analogous.
\bp\label{P:wppi}For any $j=0,\ldots,\ipartn$,
\begin{equation}  \label{E:wppi}
\pi(p_j)=
\begin{cases}
p_0 & \text{ if }j=0\,,\\
p_j - p_{j-1}& \text{ if }1<j <\ipartn\,,\\
-p_{\ipartn -1} & \text{ if }j = \ipartn\,.
\end{cases}
\end{equation}
Therefore, $\pi$ maps $\wp_n$ onto $\wp_{n-2}$ with a $1$-dimensional kernel spanned by

\[\sum_{j=0}^{\ipartn} p_j=\sum_{u\in\calSn}u\,.\]
\ep
\bpf We have
\begin{align*}
 \pi (p_j) &= \pi(\sum_{\#F=j}P_F)
= \pi\Bigl(\sumsub{1, 2 \notin
F\\ \#F=j}P_F + \sumsub{(2 \notin F)\\ \#F=j-1} P_{F \cup \{1\}} +
\sumsub{2\in F \\ \#F=j} P_F\Bigr)\\
&\equal{\eqref{E:defpi}}\sumsub{1, 2 \notin
F\\ \#F=j}P_{F-2} - \sumsub{1,2 \notin F\\ \#F=j-1} P_{F-2}=
p_j - p_{j-1}.
\end{align*}
The remaining assertions follow at once.
\epf

It follows from~\eqref{E:solbeta} that the kernel of $\beta^2:\sB\to\sol{B_{n-2}}$
 has dimension 2; in fact
\[\ker(\beta^2)=\Span\{\sum_{j=0}^n y_j,\ \sum_{j=0}^n j\cdot y_j\}=\Span
\{x_n,x_{n-1}\}\,.\]
$\varphi$ maps this kernel to the 1-dimensional space spanned by $\sum_{j=0}^{\ipartn} p_j$, since
\[\varphi(\sum_{j=0}^n y_j)=2^n\cdot\sum_{j=0}^{\ipartn} p_j
 \text{ \ \ and \ \ }
\varphi(\sum_{j=0}^n j\cdot y_j)=n2^{n-1}\cdot\sum_{j=0}^{\ipartn} p_j\,.\]
These follow from~\eqref{E:sBtowp-d} together with $\sum_{j=0}^n\binom{n}{j}=2^n$ and
$\sum_{j=0}^n j\binom{n}{j}=n2^{n-1}$.
\bc\label{C:sBexactsequence}
There is a commutative diagram
\[\xymatrix{   0\ar[r]& {\Span\{x_n,\ x_{n-1}\}}
\xyinc[r]\ar[d]_{\varphi}  &
{\sB}\ar[r]^{\beta^2}\ar[d]^{\varphi} &
{\sol{B_{n-2}}}\ar[d]^{\varphi}\ar[r] & 0\\
0\ar[r]& {\Span\{\sum_{j=0}^{\ipartn} p_j\}}
\xyinc[r] &
{\wp_n}\ar[r]_{\pi}&{\wp_{n-2}}\ar[r] & 0 }\]
where the rows are exact and the vertical maps are surjective.
\ec
\bpf
This follows from Corollary~\ref{C:Bexactsequence}, Propositions~\ref{P:sBtowp}, \ref{P:solbeta}
and~\ref{P:wppi}
and the previous remarks.
\epf

\subsection{Commutative subalgebras of the descent algebra of type D}\label{S:comm-D}

Grouping permutations of type D according to the number of descents does not
yield a subalgebra of $\SolD$, as one may easily see. On the other hand, one may
obtain commutative subalgebras of $\SolD$ simply as the image
of the commutative subalgebras of $\SolB$ under the map $\chi$ of Section~\ref{S:morphisms}.
This leads to the subalgebras considered by Mahajan in~\cite[Section 5.2]{Mah01}.
We describe the result next.

Let $\sD$, $i^{1',1}_n$ and $\hsD$ be the images under $\chi$ of $\sB$, $i^0_n$ and
$\hsB$, respectively. Thus, $\hsD$ is a commutative semisimple subalgebra of $\SolD$,
$\sD$ is a (commutative) semisimple subalgebra of $\hsD$, $i^{1',1}_n$ is a two-sided ideal
of $\hsD$ and
$\hsD=\sD+ i^{1',1}_n$. In addition, Proposition~\ref{P:BDmap} implies that
\[
\chi(x^0_j)
=\sumsub{J\subseteq [2,n-1]\\\#J=j-1}\!\!\!X_{\{1'\}\cup J}+\sumsub{J\subseteq [2,n-1]\\\#J=j-1}
\!\!\!X_{\{1\}\cup J}+ 2\cdot\!\!\!\!\!\sumsub{J\subseteq [2,n-1]\\\#J=j-2}
\!\!\!X_{\{1',1\}\cup J}=\sumsub{1'\in J\subseteq\Dsets\\\#J=j}\!\!\!X_J+\sumsub{1\in J\subseteq\Dsets\\\#J=j}
\!\!\!X_J\,,
\]
for each $j=1,\ldots,n$, and
\[
\chi(x_j)=\chi(x^0_j)+\sumsub{J\subseteq [2,n-1]\\\#J=j}\!\!\!X_J+
\sumsub{J\subseteq [2,n-2]\\\#J=j-1}\!\!\!X_{\{1',1\}\cup J}\,,
\]
for each $j=0,\ldots,n$. For instance, for $n=3$ we have
\[\chi(x^0_1)=X_{\{1'\}}+X_{\{1\}},\ \ \chi(x^0_2)=X_{\{1',2\}}+X_{\{1,2\}}+2\cdot X_{\{1',1\}},\ \
\chi(x^0_3)=2\cdot X_{\{1',1,2\}}\]
and
\[\chi(x_0)=X_\emptyset,\ \ \chi(x_1)=\chi(x^0_1)+X_{\{2\}}+X_{\{1',1\}},\ \
\chi(x_2)=\chi(x^0_2)+X_{\{1',1,2\}},\ \  \chi(x_3)=\chi(x^0_3)\,.\]

It is easy to see that both sets $\{\chi(x_j)\}_{j=0,\ldots,n}$
and $\{\chi(x^0_j)\}_{j=1,\ldots,n}$ are linearly independent and that the only linear
relation among these elements is $\chi(x_n)=\chi(x^0_n)$. Therefore, $\chi$ induces
isomorphisms
\[\sB\cong\sD,\ \  i^0_n\cong i^{1',1}_n \text{ \ and \ }\hsB\cong\hsD\,. \]
It follows from Corollary~\ref{C:BDA}, diagram~\eqref{E:gammabeta} and Corollary~\ref{C:sBexactsequence} that we have
 commutative diagrams
 \[\smallxymatrix{{\sB}\ar[rr]^{\chi}\ar[dr]_{\varphi} & &{\sD}\ar[ld]^{\psi}\\
& {\wp_n} }  \ \ \
\smallxymatrix{{\sB}\ar[rr]^{\chi}\ar[dr]_{\beta^2} & &{\sD}\ar[ld]^{\gamma}\\
& {\sol{B_{n-2}}} }\]
 and
\[\xymatrix{0\ar[r]& {\Span\{\chi(x_n),\ \chi(x_{n-1})\}}\xyinc[r]\ar@{->>}[d]_{\psi}  &
{\sD}\ar[r]^{\gamma}\ar@{->>}[d]^{\psi} & {\sol{B_{n-2}}}\ar@{->>}[d]^{\varphi}\ar[r] & 0\\
0\ar[r]& {\Span\{\sum_{j=0}^{\ipartn} p_j\}}\xyinc[r] & {\wp_n}\ar[r]_{\pi}&{\wp_{n-2}}\ar[r] & 0 }\]
where $\gamma$ is the map of~\eqref{E:defgamma}.

\section{The algebra of Mantaci and Reutenauer. Right ideals}\label{S:ManReu}

\subsection{Bases of the algebra of Mantaci and Reutenauer}\label{S:bases}
In order to recall the construction of Mantaci and Reutenauer, we need to introduce
some notation.

An (ordinary) {\em composition} of $n$ is a sequence $(a_1,\ldots,a_k)$ of positive integers such
that $a_1+\cdots +a_k=n$.
We make use of the standard bijection between compositions of $n$ and subsets of $\Asets$
\[(a_1,\ldots,a_k)\mapsto \{a_1,a_1+a_2,\ldots,a_1+\cdots+a_{k-1}\}\]
to identify the poset of subsets under inclusion with the poset of compositions
under refinement.

A {\em pseudo composition} of $n$ is a sequence $(a_0,a_1,\ldots,a_k)$ of integers such
that $a_0\geq 0$, $a_i>0$ for $i\geq 1$, and $a_0+a_1+\cdots +a_k=n$.
Pseudo composition of $n$   are in bijection
with subsets of $\Bsets$ by means of the same construction as above, and inclusion of
subsets corresponds to refinement of compositions.

Below, we will use these identifications to index basis elements of $\SolA$ or $\SolB$ by
compositions instead of subsets. Thus, for instance, if $J\subseteq\Asets$ is the subset corresponding
to the ordinary composition $\alpha=(a_1,\ldots,a_k)$, we will use $X_\alpha$ to denote
the element $X_J\in\SolA$ and $X^0_\alpha$ or $X_{(0,a_1,\ldots,a_k)}$
to denote the element $X^0_J=X_{\{0\}\cup J}\in I_n^0$.

A {\em signed composition} of $n$ is a sequence $(a_1,\ldots,a_k)$ of non-zero integers such that
$\abs{a_1}+\cdots +\abs{a_k}=n$.
There are $2\cdot 3^{n-1}$ signed compositions of $n$.

The {\em parts} of an (ordinary, pseudo or signed) composition are the integers $a_i$.

A {\em segment} of a signed composition $\alpha=(a_1,\ldots,a_k)$ is a signed composition
of the form $(a_i,a_{i+1},\ldots,a_j)$ where $[i,j]$ is a maximal interval of $[k]$ on which
the sign of the parts of $\alpha$ is constant. For instance, the segments of
$(\bar{2},1,\bar{1},\bar{2},2,2,3)$ are
\[(\bar{2})\,,\  (1)\,, (\bar{1},\bar{2})\,,\ (2,2,3)\,.\]

Following Mantaci and Reutenauer, we make the following:
\bd\label{D:OmeBbases}
Let $\alpha=(a_1,\ldots,a_k)$ be  a signed composition of $n$.
Consider the successive subintervals of $[n]$ of lengths $\abs{a_1},\ldots,\abs{a_k}$.
Define elements $\Stilde_\alpha$, $S_\alpha$ and $T_\alpha$ of the group algebra of $\calBn$ as follows:
\be
\item \label{defStilde}$\Stilde_\alpha$ is the sum of all signed permutations $w$
such that on each of those intervals the entries $w_j$ are increasing and have constant sign
equal to the sign of the corresponding part of $\alpha$.
\item \label{defS} $S_\alpha$ is the sum of all signed permutations $w$
such that on each of those intervals the numbers $\abs{w_j}$ are increasing  and have
constant sign equal to the sign of the corresponding part of $\alpha$.
\item \label{defT}
$T_\alpha$ is the sum of all signed permutations $w$ satisfying the conditions in~\eqref{defS}
and the extra requirement that those intervals are maximal with those two properties.
\ee
\ed

For instance, if $\alpha=(2,2,\bar{3},\bar{1},1)$ then these elements are
the sum of all signed permutations of the form $w=w_1w_2w_3w_4\Bar{w}_5\Bar{w}_6\Bar{w}_7\Bar{w}_8w_9$
and such that, respectively,
\begin{align*}
\Stilde_\alpha &: w_1<w_2,\ w_3<w_4,\ \Bar{w}_5<\Bar{w}_6<\Bar{w}_7\\
S_\alpha &: w_1<w_2,\ w_3<w_4,\ w_5<w_6<w_7\\
T_\alpha &: w_1<w_2>w_3<w_4,\ w_5<w_6<w_7>w_8\,.\\
\end{align*}

We will at times identify these elements with the classes of permutations that define them
(and we will do the same for the elements $X_\alpha$ of $\SolA$ or $\SolB$).

\bd\label{D:OmeB} The algebra of Mantaci and Reutenauer is the subspace $\OmeB$
 of the group algebra of $\calBn$ linearly spanned by the
elements $T_\alpha$, as $\alpha$ runs over the set of signed compositions of $n$.
\ed

\bp\label{P:OmeB} $\OmeB$ is a subalgebra of the group algebra of $\calBn$ containing
the descent algebra $\SolB$.
\ep
\bpf This is Theorem 3.9 in~\cite{ManReu}.
\epf

Given a signed permutation $w$ in the class $T_\alpha$, one may recover $\alpha$ as follows:
$\abs{a_i}$ is the length
of the $i$-th maximal interval of $[n]$ on which $\abs{w_j}$ is increasing and $\sign(w_j)$ is
constant; $\sign(a_i)$ is the sign of any $w_j$ in this interval. For instance, if $w=
\bar{3}4617\bar{5}\bar{2}\bar{8}$ then $\alpha=(\bar{1},2,2,\bar{1},\bar{2})$. Therefore,
these classes  are disjoint and the elements $T_\alpha$ form a basis of $\OmeB$;
in particular, $\dim\OmeB=2\cdot 3^{n-1}$.

The elements $S_\alpha$ also form a basis of $\OmeB$. In order to see this,
partially order the set of signed compositions of $n$ as follows:  say that $\alpha\leq\beta$ if
they have the same number of segments, these segments have the same signs,
 and the $i$-th segment of $\beta$ refines the $i$-th
segment of $\alpha$, for every $i$. For instance,
\[(\bar{2},1,\bar{3},2,5)\leq(\bar{2},1,\bar{1},\bar{2},2,2,3)\,.\]
It follows immediately from Definition~\ref{D:OmeBbases} that
\begin{equation*}
S_\alpha=\sumsub{\beta\leq\alpha}T_\beta\,.
\end{equation*}
Thus, as $\alpha$ runs over the signed compositions of $n$, the elements $S_\alpha$
form another basis of $\OmeB$.

The elements $\Stilde_\alpha$ do too. This fact is not given in~\cite{ManReu} (in fact, it is
not evident that these elements belong to $\OmeB$),
but it may be deduced as follows.

Given an ordinary composition $\alpha$ of $n$, let $\alpha^c$ denote the composition of $n$
that corresponds to the complement in $[n-1]$ of the subset corresponding to $\alpha$.
For instance, if $\alpha=(2,1,2)$ then $\alpha^c=(1,3,1)$. Given a signed composition of $n$
$\alpha$, let $\widetilde{\alpha}$ denote the signed composition obtained by replacing each negative
segment $\alpha_i$ of $\alpha$ by $\alpha_i^c$ (viewing it as an ordinary composition)
and keeping its (negative) sign. For instance,
if $\alpha=(3,\underbrace{\bar{2},\bar{1},\bar{2}}_{\alpha_1},4,2,\underbrace{\bar{3}}_{\alpha_2},1)$ then
$\widetilde\alpha=(3,\underbrace{\bar{1},\bar{3},\bar{1}}_{\alpha_1^c},4,2,
\underbrace{\bar{1},\bar{1},\bar{1}}_{\alpha_2^c},1)$. Note that $\alpha\mapsto\widetilde\alpha$ is
an involution on the set of signed compositions of $n$.

Introduce a new partial order on the set of signed compositions of $n$ by declaring that
$\alpha\preceq\beta$ if
they have the same number of segments, these segments have the same signs,
and the $i$-th segment of $\beta$ refines the $i$-th
segment of $\alpha$ if they are both positive and  the $i$-th segment of $\alpha$ refines the $i$-th
segment of $\beta$ if they are both negative, for every $i$. For instance,
\[(\bar{2},1,\bar{1},\bar{2},2,2,3)\preceq(\bar{2},1,\bar{3},2,1,1,3)\,.\]
Now, it is easy to see that with the above definitions
\[\Stilde_\alpha=\sumsub{\beta\preceq\widetilde{\alpha}}T_\beta\,.\]
This immediately implies that the elements $\Stilde_\alpha$ form a basis of $\OmeB$.

We summarize these facts in the following lemma.
\bl\label{L:OmeBbases}
As $\alpha$ runs over the signed compositions of $n$, the sets $\{T_\alpha\}$,
$\{S_\alpha\}$ and $\{\Stilde_\alpha\}$ are
 bases of $\OmeB$.
\el

\subsection{The image of the algebra of Mantaci and Reutenauer}\label{S:OmetoSol}

We now consider the map $\varphi:\calBn\to\calSn$ that forgets the signs,
restricted to $\OmeB$. To describe its action on the various bases of $\OmeB$,
we need to introduce some notation.

Given a signed composition $\alpha=(a_1,\ldots,a_k)$ of $n$, we construct four
ordinary compositions of $n$.
First, $\abs{\alpha}:=(\abs{a_1},\ldots,\abs{a_k})$ denotes the ordinary composition of $n$
obtained by forgetting the signs of the parts of $\alpha$. Another ordinary
composition $\underline{\alpha}$ is constructed as follows.  Split $[k]$ into maximal intervals
$[i,j]$ with the property that the signs of the parts $a_i,a_{i+1},\ldots,a_j$ alternate.
The parts of $\underline{\alpha}$ are then $\abs{a_i}+\abs{a_{i+1}}+\ldots+\abs{a_j}$, one for each
such interval, in the same order as the intervals appear. For instance, if
$\alpha=(\underbrace{\bar{2},1,\bar{1}},\underbrace{\bar{2},2},\underbrace{2},\underbrace{3})$
then $\underline{\alpha}=(4,4,2,3)$. Two other ordinary compositions $\Oalpha$ and
$\Ualpha$  are obtained by the following procedures. Consider the various segments
$\alpha_i$ of $\alpha$ (Section~\ref{S:bases}) and let $n_i$ be the sum of the parts of $\alpha_i$.
 To obtain $\Oalpha$, simply replace
each negative segment $\alpha_i$ by $(1,\ldots,1)$ (with $n_i$ parts equal to $1$)
and drop the signs.
To obtain $\Ualpha$, replace each
negative segment $\alpha_i$ by its complement $\alpha_i^c$ (keeping the negative signs),
each positive segment $\alpha_j$ by $(n_j)$, and then apply the operation
$\alpha\mapsto\underline{\alpha}$.
For instance, if $\alpha=
(\bar{2},1,\bar{1},\bar{2},2,2,3)$ then
$\Oalpha=(1,1,1,1,1,1,2,2,3)$ and $\Ualpha=(1,4,8)$.

With this notation, we have:

\bp\label{P:OmetoSol} For any signed composition $\alpha$ we have
\begin{equation}\label{E:OmetoSol}
\varphi(S_\alpha)=X_{\abs{\alpha}}\,,\ \
\varphi(T_\alpha)=\sum_{\underline{\alpha}\leq\beta\leq\abs{\alpha}}Y_\beta
\text{ \ and \ }\varphi(\Stilde_\alpha)=\sum_{\Ualpha\leq\beta\leq\Oalpha}Y_\beta\,.
\end{equation}
In particular, $\varphi$ maps $\OmeB$ onto $\SolA$.
\ep
\bpf Write $\alpha=(a_1,\ldots,a_k)$ and
consider the successive subintervals of $[n]$ of lengths $\abs{a_1},\ldots,\abs{a_k}$ as in
Definition~\ref{D:OmeBbases}. According to the definition of descents for ordinary permutations,
$X_{\abs{\alpha}}$ is the sum of those permutations $u\in\calSn$ such that the $u_j$
are increasing on each such interval. Clearly, for each such $u$
there is a unique $w\in\calBn$ satisfying condition~\eqref{defS} in Definition~\ref{D:OmeBbases} and
 $\varphi(w)=u$. Hence,
$\varphi(S_\alpha)=X_{\abs{\alpha}}$. This implies that $\varphi(\OmeB)=\SolA$. The
expressions for $\varphi(T_\alpha)$ and $\varphi(\Stilde_\alpha)$ follow
from a similar analysis that we omit.
\epf

The previous result is of crucial importance. It entails that, from the point of view of
this paper, the ``correct'' analog of the descent algebra of type A is the algebra $\OmeB$ of
Mantaci and Reutenauer. On the other hand, the descent algebra of type B is, one may say,
the type B analog of the peak algebra. The precise meaning of these assertions is
summarized by the commutative diagram
\[\xymatrix@C=0pc{ {I_n^{0}}\ar@{->>}[d]_{\varphi} &\subseteq &
{\SolB}\ar@{->>}[d]_{\varphi} &\subseteq& {\OmeB}\ar@{->>}[d]^{\varphi} &\subseteq&
 {\Q\calBn}\ar@{->>}[d]^{\varphi}\\
 {\pppnint} &\subseteq& {\pppn}& \subseteq & {\SolA} &\subseteq& {\Q\calSn} }\]

\subsection{The canonical ideal in type B and the peak ideal as principal right ideals}
\label{S:principal}

Consider the standard embedding
\[\calB{p}\times\calS{q}\inc\calB{p+q},\ (w\times u)(i)=\begin{cases} w(i)
 & \text{ if }1\leq i\leq p\,,\\
 p+u(i-p) & \text{ if }p+1\leq i\leq p+q\,. \end{cases}\]
Let $\alpha=(a_1,\ldots,a_k)$ be an ordinary
composition of $n$. The obvious extension of the above embedding allows us to view
$\calB{0}\times \calS{a_1}\times\cdots\times\calS{a_k}$ as a subgroup of $\calBn$.
Consider the successive subintervals of $[n]$ of lengths
$a_1,\ldots,a_k$.
 Recall that $X^0_\alpha$ denotes the
element $X_{(0,a_1,\ldots,a_k)}\in I_n^0$. This is the sum of all signed permutations
$\xi\in\calBn$ such that, on each of those intervals, the entries $\xi_j$ are increasing.
These permutations are the {\em $(0,a_1,\ldots,a_k)$-shuffles of type B}. It is well-known,
and easy to see,
that this is a set of left representatives for the cosets of
$\calB{0}\times \calS{a_1}\times\cdots\times\calS{a_k}$ as a subgroup of $\calBn$.

In particular, the element $X^0_{(n)}=X_{(0,n)}\in I_n^0$ is the sum of all signed permutations
$\xi\in\calBn$ such that
\begin{equation}\label{E:0nshuffles}
 \xi_1<\xi_2<\cdots<\xi_n\,.
\end{equation}

The basis $\Stilde_\alpha$ of $\OmeB$ is particularly useful for our purposes
because of the following simple multiplication formula with $X^0_{(n)}$.

\bp\label{P:BStilde} For any signed composition $\alpha$ of $n$ we have
\begin{equation}\label{E:BStilde}
X^0_{(n)}\cdot\Stilde_{\alpha}=X^0_{\abs{\alpha}}\,.
\end{equation}
\ep
\bpf It is possible to deduce this result from the multiplication rule for the basis
$\Stilde_\alpha$ given in~\cite[Corollary 5.3]{ManReu}, but it seems simpler to provide a
direct argument.

Write $\alpha=(a_1,a_2,\ldots,a_k)$. Consider the three classes of signed permutations which define
the elements $X^0_n$, $\Stilde_\alpha$ and $X^0_{\abs{\alpha}}$.  If $\xi$ is in the
class $X^0_n$ and $w$ is in the class $\Stilde_\alpha$, then $\xi$ is increasing on $[n]$
and $w$ is (in particular) increasing on each of the successive subintervals of $[n]$ of lengths
$\abs{a_1},\ldots,\abs{a_k}$, according to~\eqref{E:0nshuffles} and Definition~\ref{D:OmeBbases}.
Therefore, $\xi\cdot w$ is also increasing on these intervals, i.e., $\xi\cdot w$ is a permutation in the
class $X^0_{\abs{\alpha}}$. Thus, to prove the result it suffices to check that all permutations
in $X^0_{(n)}\cdot\Stilde_{\alpha}$ are distinct and that the total
number agrees with the number of permutations in the class $X^0_{\abs{\alpha}}$.

According to the previous remarks, the number of permutations in $X^0_{(n)}$ and $X^0_{\abs{\alpha}}$
are respectively
\[\frac{2^n\cdot n!}{n!} \text{ \  and \ }\frac{2^n\cdot n!}{\abs{a_1}!\cdots\abs{a_k}!}\,.\]
On the other hand, it follows from Definition~\ref{D:OmeBbases} that, after reversing the
order of the entries within each negative segment and dropping the signs,
the permutations $w$ which
appear in the class $\Stilde_{\alpha}$ become precisely the permutations $u$ in the
class $X_{\abs{\alpha}}\in\SolA$. This is the class of $(\abs{a_1},\ldots,\abs{a_k})$-shuffles of
type A, which is a set of representatives for the subgroup
$\calS{\abs{a_1}}\times\cdots\times\calS{\abs{a_k}}$ of $\calSn$. Hence,  the number of
permutations in this class is
\[ \frac{n!}{\abs{a_1}!\cdots\abs{a_k}!}\,.\]
Comparing these quantities we see that the total number of permutations on both sides
of~\eqref{E:BStilde} is the same.

It only remains to check that repetitions do not occur in the product
$X^0_{(n)}\cdot\Stilde_{\alpha}$. Suppose $\xi\cdot w=\xi'\cdot w'$ with
$\xi,\xi'\in X^0_{(n)}$ and $w,w'\in\Stilde_{\alpha}$. Applying the morphism
 $\varphi:\calBn\to\calSn$ we deduce that for each $i=1,\ldots,n$,
 \begin{equation}\tag{A}
 \abs{\xi_{\abs{w_i}}}=\abs{\xi'_{\abs{w'_i}}}\,.
 \end{equation}
On the other hand,
\[\sign\bigl((\xi\cdot w)_i\bigr)=\sign(\xi_{\abs{w_i}})\sign(w_i)
\text{ \ and \ } \sign\bigl((\xi'\cdot w')_i\bigr)=\sign(\xi'_{\abs{w'_i}})\sign(w'_i)\]
(this holds for any pair of signed permutations). Now, since $w,w'\in\Stilde_{\alpha}$,
 Definition~\ref{D:OmeBbases} implies that $\sign(w_i)=\sign(w'_i)$. Hence,
 \begin{equation}\tag{B}
 \sign(\xi_{\abs{w_i}})=\sign(\xi'_{\abs{w'_i}})\,.
 \end{equation}
Equations (A) and (B) imply that $\xi_{\abs{w_i}}=\xi'_{\abs{w'_i}}$ for every $i$.
But $\xi$ and $\xi'$ are
representatives for $\calSn$ in $\calBn$, so we must have $\xi=\xi'$, and then
$w=w'$. This completes the proof.
\epf

\bc\label{C:rightOmeB}  $I_n^0$ is a right ideal of $\OmeB$. Moreover, it is the
principal right ideal generated by $X^0_{(n)}$.
\ec
\bpf As $\alpha$ runs over the signed compositions of $n$,
 the elements $\Stilde_\alpha$ form a basis of $\OmeB$
(Lemma~\ref{L:OmeBbases})
and the elements $X^0_{\abs{\alpha}}$ span $I_n^0$. Thus, by Proposition~\ref{P:BStilde},
\[X^0_{(n)}\cdot\OmeB=I_n^0\,.\]
\epf

 Consider now the element $\Pint_{\emptyset}\in\pppnint$. In order to make explicit the dependence
on $n$, we will denote it by $\Pint_{(n)}$. According to~\eqref{E:defpeakintone}, $\Pint_{(n)}$
is the sum of all permutations $u\in\calSn$ such that
\begin{equation}\label{E:nointpeaks}
u_1>\cdots>u_i<u_{i+1}<\cdots<u_n \text{ \ for some \ }i=1,\ldots,n\,.
\end{equation}
Comparing~\eqref{E:0nshuffles} with~\eqref{E:nointpeaks} (or using~\eqref{E:Btointpeaks})
we see that
\begin{equation}\label{E:shufflesnopeaks}
\varphi(X^0_{(n)})=2\cdot\Pint_{(n)}\,.
\end{equation}

\bc\label{C:rightSolA}$\pppnint$ is a right ideal of $\SolA$. Moreover, it is the
principal right ideal generated by $\Pint_{(n)}$.
\ec
\bpf By Corollary~\ref{C:rightOmeB}, $I_n^0=X^0_{(n)}\cdot\OmeB$. Applying $\varphi$
and using Theorem~\ref{T:Bideals}, Proposition~\ref{P:OmetoSol} and~\eqref{E:shufflesnopeaks} we deduce that
\[\pppnint=\Pint_{(n)}\cdot\SolA\,.\]
\epf

\br \label{R:rightideal}
\be
\item $\pppnint$ is {\em not} a left ideal of
$\SolA$. This occurs already for $n=3$, as one may easily verify that
$Y_{\{1\}}\cdot\Pint_{\{2\}}\notin\pppint{3}$. It follows that $I_n^0$ is not a left ideal
of $\OmeB$.
\item Corollary~\ref{C:rightSolA} was first obtained by Schocker~\cite[Theorems 1 and 5]{Sch}, by
other means. (Schocker endows the group algebra with the opposite product of the composition of permutations.
His left ideals are therefore right ideals in our sense.)
\item We will show below that $I_n^0$ is still generated by $X^0_{(n)}$ as a right
ideal of the smaller algebra $\SolB$, and similarly that $\pppnint$ is
generated by $\Pint_{(n)}$ as a right ideal of $\pppn$ (Proposition~\ref{P:rightpeak}).
\item We will say more about the elements $X^0_{(n)}$ and $\Pint_{(n)}$ in Section~\ref{S:words}.
\ee
\er

\subsection{The descents-to-peaks transform and its type B analog}\label{S:thetamaps}

\bd\label{D:Theta} Let $\Theta:\SolA\to\SolA$ be the unique morphism of right $\SolA$-modules
such that
\begin{equation}\label{E:deftheta}
\Theta(1)=2\cdot\Pint_{(n)}\,.
\end{equation}
Let $\ThetaB:\OmeB\to\OmeB$ be the unique morphism of right $\OmeB$-modules
such that
\begin{equation}\label{E:defthetaB}
\ThetaB(1)= X^0_{(n)}\,.
\end{equation}
We refer to $\Theta$ as the {\em descents-to-peaks transform} and to $\ThetaB$ as its type B
analog.
\ed

\br\label{R:thetadef} The descents-to-peaks transform should be attributed to
Stembridge, and independently, to Krob, Leclerc and Thibon, as we now explain.
\be
\item Below~\eqref{E:theta} we will derive an expression for $\Theta$ which implies
that this map is dual to a map previously considered by Stembridge~\cite{Ste}
(the precise connection  will be worked out in detail in Section~\ref{S:duality}).
\item Krob, Leclerc and Thibon have considered a family of maps $\SolA\to\SolA$ which
depend on a parameter $q$~\cite[Section 5.4]{KLT}. Comparing our Definition~\ref{D:Theta}
with~\cite[Proposition 5.2]{KLT} we see that our map $\Theta$ is the specialization of
their map at $q=-1$. Many interesting results about $\Theta$ are given in this reference.
In particular, it is easy to see that (the case $q=-1$ of)~\cite[Proposition 5.41]{KLT}
is equivalent to our~\eqref{E:theta}.
Our work adds the observation on the connection with the map of Stembridge, as well as
the introduction of the type B analog.
\ee
\er

It follows formally from the definitions, together with~\eqref{E:shufflesnopeaks} and the fact that
$\varphi$ is a surjective morphism of algebras, that the diagram
\begin{equation}\label{E:thetadiagram}
\xymatrix{ {\OmeB}\ar[r]^{\ThetaB}\ar[d]_{\varphi} & {\OmeB}\ar[d]_{\varphi}\\
 {\SolA}\ar[r]_{\Theta} & {\SolA} }
\end{equation}
commutes.

We may easily deduce the following explicit formulas for these maps.

\bp\label{P:theta} For any signed composition $\alpha$ of $n$ we have
\begin{equation}\label{E:thetaB}
\ThetaB(\Stilde_\alpha)=X^0_{\abs{\alpha}}\,,
\end{equation}
and for any subset $J\subseteq\Asets$ we have
\begin{equation}\label{E:theta}
\Theta(X_J)=2^{1+\#J}\cdot\!\!\!\!\sumsub{F\in\calFnint\\F\subseteq J\cup(J+1)}\Pint_{F}\,.
\end{equation}
In particular, the image of $\ThetaB$ is $I_n^0$ and the image
of $\Theta$ is $\pppnint$.
\ep
\bpf Since $\ThetaB$ is a map of right $\OmeB$-modules,
\[\ThetaB(\Stilde_\alpha)=\ThetaB(1)\cdot \Stilde_\alpha
\equal{\eqref{E:defthetaB}} X^0_{(n)}\cdot \Stilde_\alpha\equal{\eqref{E:BStilde}}
X^0_{\abs{\alpha}}\,.\]
Let $\alpha$ be the ordinary composition of $n$ corresponding to the subset $J$. We may
also view it as a signed composition of $n$ (with positive signs). For such compositions
we have $\Stilde_\alpha=S_\alpha$, according to Definition~\ref{D:OmeBbases}.
Therefore, by ~\eqref{E:OmetoSol}, $\varphi(\Stilde_\alpha)=X_J$. Now, from the
commutativity of~\eqref{E:thetadiagram} we deduce
\[\Theta(X_J)=\Theta\varphi(\Stilde_\alpha)=\varphi\ThetaB(\Stilde_\alpha)
\equal{\eqref{E:thetaB}}\varphi(X^0_J)\equal{\eqref{E:Btointpeaks}}
2^{1+\#J}\cdot\!\!\!\!\sumsub{F\in\calFnint\\F\subseteq J\cup(J+1)}\Pint_{F}\,,\]
as claimed in~\eqref{E:theta}.

Since $\ThetaB$ maps a basis of $\OmeB$  onto a basis of $I_n^0$,
the assertion about its image follows. The image of $\Theta$ is therefore
$\varphi(I_n^0)=\pppnint$ (Theorem~\ref{T:Bideals}).
\epf

\bp\label{P:thetabijective} The restriction of $\ThetaB$ to $I_n^0$ is a bijective map
$I_n^0\to I_n^0$. The restriction of $\Theta$ to $\pppnint$ is a bijective map
$\pppnint\to\pppnint$.
\ep
\bpf Let $\alpha=(a_1,\ldots,a_k)$ be an ordinary composition of $n$ and consider the
class of signed permutations $X^0_\alpha$. A permutation $w$ in this class is increasing
on each of the successive subintervals of $[n]$ of lengths $a_1,\ldots,a_k$. Let $[t_i,t_i+a_i]$
be one such interval. There is an index $b_i=0,\ldots,a_i$ such that
\[w_{t_i}<\ldots<w_{t_i+b_i}<0<w_{t_i+b_i+1}<\ldots<w_{t_i+a_i}\,.\]
Consider the signed composition
\[\alpha_b:=(\bar{b}_1,a_1-b_1,\bar{b}_2,a_2-b_2,\ldots,\bar{b}_k,a_k-b_k)\]
where we understand that any possible zero parts are omitted.
Definition~\ref{D:OmeBbases} implies that $w$ belongs to the class $\Stilde_{\alpha_b}$.
This shows that the class $X^0_\alpha$ is the union of the
classes $\Stilde_{\alpha_b}$, over all possible choices of indeces $b_i=0,\ldots,a_i$.
This union is disjoint, since permutations in distinct classes $\Stilde_{\alpha_b}$
are distinguished by the signs of their entries.
Thus, we may write
\[ X^0_\alpha=\sum_{b_i=0}^{a_i}\Stilde_{\alpha_b}\,.\]
Now, by~\eqref{E:thetaB},
\[\ThetaB(X^0_\alpha)=\sum_{b_i=0}^{a_i} X^0_{\abs{\alpha_b}}\,.\]
The extreme choices $b_i=0$ or $a_i$ $\forall i$ yield $\abs{\alpha_b}=\alpha$.
These are $2^k$ choices. All other choices yield compositions $\alpha_b$ that are
finer than $\alpha$. Therefore,
\[\ThetaB(X^0_\alpha)=2^k\cdot X^0_\alpha +\sum_{\alpha<\beta}c_\beta\cdot X^0_\beta\]
for some non-negative integers $c_\beta$. This implies that $\ThetaB:I_n^0\to I_n^0$ is
bijective.

The assertion about $\Theta$ follows from the assertion about $\ThetaB$ and the
surjectivity of $\varphi$, as usual.
\epf
 We can now deduce that the canonical ideal in type B and the peak ideal
are principal ideals of $\SolB$ and $\pppn$. In fact, they are principal
as right ideals. Note that these properties are neither stronger nor weaker than
those in Corollaries~\ref{C:rightOmeB} and~\ref{C:rightSolA}.

\bp\label{P:rightpeak} As a right ideal of $\SolB$, $I_n^0$ is generated by $X^0_{(n)}$.
As a right ideal of $\pppn$, $\pppnint$ is generated by $\Pint_{(n)}$.
\ep
\bpf Since $I_n^0\subseteq\SolB\subseteq\OmeB$, Propositions~\ref{P:thetabijective}
and~\ref{P:theta} imply that
\[\ThetaB(\SolB)=I_n^0\,.\]
 Hence,
\[I_n^0=\ThetaB(\SolB)=\ThetaB(1)\cdot\SolB= X^0_{(n)}\cdot\SolB\,.\]
The assertion about $\pppnint$ follows similarly or by applying $\varphi$.
\epf

\begin{rem} It is also possible to deduce the result about $I_n^0$ from the multiplication rule for the
$X$-basis of $\SolB$ given in~\cite[Theorem 1]{BB92a}.
\end{rem}

\section{Hopf algebraic structures}\label{S:external}

In this section we consider the direct sum over all $n\geq 0$ of the group, descent and
peak algebras discussed elsewhere in the paper. These spaces may be endowed with a new
product (called the external product) and a coproduct, so that in the diagram
\[\xymatrix@C=0pc{ {I^{0}}\ar@{->>}[d]_{\varphi} &\subseteq &
{\Sol{B}}\ar@{->>}[d]_{\varphi} &\subseteq& {\Ome{B}}\ar@{->>}[d]^{\varphi} &\subseteq&
 {\QB}\ar@{->>}[d]^{\varphi}\\
 {\pppint{}} &\subseteq& {\ppp{}}& \subseteq & {\Sol{A}} &\subseteq& {\QS} }\]
 all objects are Hopf algebras, except for $\Sol{B}$, which is an $I^0$-module coalgebra,
 and $\ppp{}$, which is a $\pppint{}$-module coalgebra.

\subsection{The Hopf algebras of permutations and signed permutations}\label{S:permutations}

Let $p$ and $q$ be non-negative integers and $n=p+q$.
Consider the standard embedding
\[\calS{p}\times\calS{q}\inc\calS{p+q},\ (u\times v)(i)=\begin{cases} u(i)
 & \text{ if }1\leq i\leq p\,,\\
 p+v(i-p) & \text{ if }p+1\leq i\leq p+q\,. \end{cases}\]

The set of $(p,q)$-shuffles is
\[\Sh(p,q)=\{\xi\in\calSn\ \mid\ \xi_1<\cdots<\xi_p \text{ \ and \ }\xi_{p+1}<\cdots<\xi_{n}\}\,.\]
This is a set of left representatives for the cosets of
$\calS{p}\times \calS{q}$ as a subgroup of $\calSn$. In other words, given $p$ and $w\in\calSn$,
there is a unique triple $(\xi,w_{(p)},w'_{(p)})$ such that $\xi\in\Sh(p,n-p)$, $w_{(p)}\in\calS{p}$,
 $w'_{(p)}\in\calS{n-p}$ and
 \begin{equation}\label{E:parabolic}
 w=(w_{(p)}\times w'_{(p)})\cdot\xi^{-1}\,.
 \end{equation}

Consider the direct sum of all group algebras of the symmetric groups
 \[\QS:=\Q\oplus\Q\calS{1}\oplus\Q\calS{2}\oplus\cdots\,.\]
This space carries a graded Hopf algebra structure defined as follows. Given $u\in\calS{p}$ and
$v\in\calS{q}$, their product is
\[u\ast v:=\sum_{\xi\in\Sh(p,q)}\xi\cdot(u\times v)\,.\]
The coproduct of $w\in\calSn$ is
\[\Delta(w):=\sum_{p=0}^n w_{(p)}\ten w'_{(p)}\]
where $w_{(p)}$ and $w'_{(p)}$ are as in~\eqref{E:parabolic}.

This is the Hopf algebra of Malvenuto and Reutenauer~\cite{MalReu}.
A detailed study of this Hopf algebra is given in~\cite{AS} (the object
denoted by $\SSym$ in this reference is the dual Hopf algebra to $\QS$).

Similar constructions exist for the groups of signed permutations. First, we have
the embedding $\calB{p}\times\calB{n-p}\inc\calBn$, $(u,v)\mapsto u\times v$, defined
by the same formulas as before. The set of $(p,n-p)$-shuffles of type A is still a set
of left representatives for the cosets of $\calB{p}\times\calB{n-p}$ as a subgroup of $\calBn$,
which allows us to associate signed permutations $w_{(p)}\in\calB{p}$
and $w'_{(p)}\in\calB{n-p}$ to a signed permutation $w\in\calBn$ by means of~\eqref{E:parabolic}.

The space
\[\QB:=\Q\oplus\Q\calB{1}\oplus\Q\calB{2}\oplus\cdots\,.\]
 carries a graded Hopf algebra structure with product
\[u\ast v:=\sum_{\xi\in\Sh(p,q)}\xi\cdot(u\times v)\]
(for $u\in\calB{p}$ and $v\in\calB{q}$) and coproduct
\[\Delta(w):=\sum_{p=0}^n w_{(p)}\ten w'_{(p)}\]
(for $w\in\calBn$). We refer the reader
to~\cite{AM} for a proof of this fact and for a detailed study of the Hopf algebra
 of signed permutations $\QB$.

When an explicit distinction is necessary, one refers to the product $u\ast v$ of $\QS$ (and $\QB$) as the
{\em external} product, and to the product $u\cdot v$ in the group algebras $\Q\calSn$ (and $\Q\calBn$)
as the {\em internal} product.

It is clear from the definitions that the map $\varphi:\QB\to\QS$ that forgets the
signs preserves all the structure; in particular, it is a morphism of graded Hopf algebras.

\subsection{The Hopf algebra of non-commutative symmetric functions and its type B analog}\label{S:nsym}
Consider the direct sum of all descent algebras of type A
\[\Sol{A}:=\Q\oplus\Sol{A_0}\oplus\Sol{A_1}\oplus\cdots\,.\]
As in Section~\ref{S:ManReu} we index the $X$-basis of $\SolA$ by compositions of $n$,
instead of subsets of $\Asets$.

 Since each $\SolA$ is a subspace
of $\Q\calSn$, we may view $\Sol{A}$ as a graded subspace of $\QS$.
It turns out that $\Sol{A}$ is a graded Hopf subalgebra of
$\QS$~\cite[Theorem 3.3]{MalReu}.
 Moreover, the external product of $X$-basis elements is simply
\begin{equation}\label{E:prodSolA}
X_{(a_1,\ldots,a_k)}\ast X_{(b_1,\ldots,b_h)}=X_{(a_1,\ldots,a_k,b_1,\ldots,b_h)}\,,
\end{equation}
which shows that $\Sol{A}$ is freely generated as an algebra by the elements $\{X_{(n)}\}_{n\geq 1}$,
and the coproduct is determined by
\begin{equation}\label{E:coprodSolA}
\Delta(X_{(n)})=\sum_{i+j=n}X_{(i)}\ten X_{(j)}\,,
\end{equation}
where we agree that $X_{(0)}$ denotes the unit element $1\in\Q$. Thus, $\Sol{A}$ is the Hopf
algebra of {\em non-commutative symmetric functions}~\cite{GKal,KLT}.

While the coproduct $\Delta$ is not compatible with the internal product of $\QS$
(in general, $\Delta(u\cdot v)\neq\Delta(u)\cdot\Delta(v)$), it is known that it is compatible
with the internal product of $\Sol{A}$~\cite[Remarque 5.15]{Mal}.

Similarly, the subspace
\[\Ome{B}:=\Q\oplus\Ome{B_1}\oplus\Ome{B_2}\oplus\cdots\,.\]
is a graded Hopf subalgebra of $\QB$. We will describe this structure on the
$\Stilde$-basis of Section~\ref{S:bases}.
The external product is again given by concatenation (of signed compositions)
\begin{equation}\label{E:prodOmeB}
\Stilde_{(a_1,\ldots,a_k)}\ast \Stilde_{(b_1,\ldots,b_h)}=\Stilde_{(a_1,\ldots,a_k,b_1,\ldots,b_h)}\,.
\end{equation}
It follows that $\Ome{B}$ is freely generated as an algebra by the elements
$\{\Stilde_{(n)},\ \Stilde_{(\Bar{n})}\}_{n\geq 1}$. For this reason, $\OmeB$ may be viewed as an
algebra of {\em non-commutative symmetric functions of type B}. The coproduct is determined by
\begin{equation}\label{E:coprodOmeB}
\Delta(\Stilde_{(n)})=\sum_{i+j=n}\Stilde_{(i)}\ten \Stilde_{(j)} \text{ \ and \ }
\Delta(\Stilde_{(\Bar{n})})=\sum_{i+j=n}\Stilde_{(\Bar{i})}\ten \Stilde_{(\Bar{j})}\,.
\end{equation}
These facts are simple consequences of the definitions; more details will be found in~\cite{AM}.

\subsection{The Hopf algebra of interior peaks and its type B analog}\label{S:pint}

Consider the spaces
\[\pppint{}:=\Q\oplus\pppint{1}\oplus\pppint{2}\oplus\cdots
\text{ \ and \ }
I^0:=\Q\oplus I^0_{1}\oplus I^0_{2}\oplus\cdots\,.\]
$\pppint{}$ is a graded subspace of $\Sol{A}$ and $I^0$ is a graded subspace of $\Ome{B}$.
Recall the basis $X^0_\alpha$ of $I^0_n$, indexed by ordinary compositions of $n$ as in
Section~\ref{S:bases}. It follows easily from the definitions that
\begin{equation}\label{E:prodI}
X^0_{(a_1,\ldots,a_k)}\ast X^0_{(b_1,\ldots,b_h)}=X^0_{(a_1,\ldots,a_k,b_1,\ldots,b_h)}\,,
\end{equation}
and
\begin{equation}\label{E:coprodI}
\Delta(X^0_{(n)})=\sum_{i+j=n}X^0_{(i)}\ten X^0_{(j)}\,,
\end{equation}
where we agree that $X^0_{(0)}=1\in\Q$.
\bp\label{P:hopfpint} The space $I^0$ is a graded Hopf subalgebra of $\Ome{B}$ and the space
$\pppint{}$ is a graded Hopf subalgebra of the Hopf algebra of non-commutative symmetric
functions $\Sol{A}$.
\ep
\bpf The preceding formulas show that $I^0$ is closed for the coproduct and the external product
of $\Ome{B}$.
The assertion for $\pppint{}$ follows  since $\varphi:\Ome{B}\to\Sol{A}$ is
a morphism of graded Hopf algebras such that $\varphi(I^0)=\pppint{}$ (Theorem~\ref{T:Bideals}).
\epf
In fact, formulas~\eqref{E:prodSolA}, \eqref{E:coprodSolA}, \eqref{E:prodI} and~\eqref{E:coprodI}
show that $I^0$ and  $\Sol{A}$ are isomorphic as graded Hopf algebras. On the other hand, the internal products
of $I_n^0$ and $\SolA$ are of course distinct, and as is clear from the results of this paper,
these two objects play very different roles and should be kept apart.

Recall the maps $\ThetaB:\OmeB\to I^0_n$ and $\Theta:\SolA\to\pppnint$ of Section~\ref{S:thetamaps}.
Summing over all $n\geq 0$ they give rise to maps $\Ome{B}\to I^0$ and $\Sol{A}\to\pppint{}$
that we still denote by the same symbols.

\bp\label{P:thetahopf} The maps  $\ThetaB:\Ome{B}\to I^0$ and $\Theta:\Sol{A}\to\pppint{}$
are morphisms of graded Hopf algebras.
\ep
\bpf By Proposition~\ref{P:theta}, $\ThetaB(\Stilde_\alpha)=X^0_{\abs{\alpha}}$. Comparing
~\eqref{E:prodSolA} and~\eqref{E:coprodSolA} with~\eqref{E:prodI} and~\eqref{E:coprodI} we
immediately see
that $\ThetaB$ is a morphism of Hopf algebras. The assertion for $\Theta$ follows from the
commutativity of~\eqref{E:thetadiagram}.
\epf

\subsection{The module coalgebra of peaks and its type B analog}
Consider, finally, the spaces
\[\ppp{}:=\Q\oplus\ppp{1}\oplus\ppp{2}\oplus\cdots
\text{ \ and \ }
\Sol{B}:=\Q\oplus \Sol{B_1}\oplus \Sol{B_2}\oplus\cdots\,.\]
$\ppp{}$ is a graded subspace of $\Sol{A}$ which contains $\pppint{}$
and $\Sol{B}$ is a graded subspace of $\Ome{B}$ which contains $I^0$.

The first observation is that $\ppp{}$ is {\em not} closed under the product of $\Sol{A}$ and $\Sol{B}$
is not closed under the product of $\Ome{B}$. For instance, the element $P_{\{1\}}\in\ppp{2}$
satisfies
\[ P_{\{1\}}\ast P_{\{1\}}=Y_{\{1,2,3\}}+Y_{\{1,3\}}\notin\ppp{4}\,.\]

However, it is easy to see that for any pseudo composition $(a_0,a_1,\ldots,a_k)$ and
any ordinary composition $(b_1,\ldots,b_h)$,
\begin{equation}\label{E:prodSolB}
X_{(a_0,a_1,\ldots,a_k)}\ast X^0_{(b_1,\ldots,b_h)}=X_{(a_0,a_1,\ldots,a_k,b_1,\ldots,b_h)}\,.
\end{equation}
Also, for any $n\geq 0$,
\begin{equation}\label{E:coprodSolB}
\Delta(X_{(n)})=\sum_{i+j=n}X_{(i)}\ten X_{(j)}\,.
\end{equation}

For the notion of module coalgebra we refer to~\cite[Definition IX.2.1]{Kas}.
\bp\label{P:hopfpeak} The space $\Sol{B}$ is a right $I^0$-module subcoalgebra of $\Ome{B}$ and
the space $\ppp{}$ is a a right $\pppint{}$-module subcoalgebra of $\Sol{A}$.
\ep
\bpf
Since $I^0$ is a Hopf subalgebra of $\Ome{B}$, one may view $\Ome{B}$ as a right $I^0$-module
coalgebra, simply by right multiplication.
Equation~\eqref{E:prodSolB} shows that $\Sol{B}$ is an $I^0$-submodule, and moreover,
that it is the free right $I^0$-module generated by the elements $\{X_{(n)}\}_{n\geq 0}$.
Together with~\eqref{E:coprodSolB}, this guarantees that $\Sol{B}$ is closed under $\Delta$
(since $\Delta$ is a morphism of algebras).
Thus, $\Sol{B}$ is a right $I^0$-module subcoalgebra of $\Ome{B}$. The assertion for $\ppp{}$
follows by applying the morphism $\varphi$, as usual.
\epf

\begin{rem} The space $\Sol{B}$, endowed with this coalgebra structure, has been considered
in the recent thesis of Chow~\cite{Chow}. Chow also considers an equivalent of the module structure
(over $\Sol{A}$ instead of $I^0$).
\end{rem}

Recall the map $\beta:\SolB\to\Sol{B_{n-1}}$ from Section~\ref{S:ideal-B}. This is a
morphism of algebras with respect to the internal structure. We use the same symbol to denote the
map $\Sol{B}\to\Sol{B}$ that vanishes on the component of degree $0$ and sends
the component $\SolB$ of degree $n\geq 1$ to $\Sol{B_{n-1}}$. Thus, $\beta:\Sol{B}\to\Sol{B}$
is a homogeneous map of degree $-1$. We now describe the
behavior of this map with respect to the external structure. Recall that any coalgebra $C$
can be viewed as a $C$-bicomodule with right and left coaction given by the coproduct (for the
definition of comodules see~\cite[Definition III.6.1]{Kas}).

\bp\label{P:beta-external} The map $\beta:\Sol{B}\to\Sol{B}$ is a morphism of $\Sol{B}$-bicomodules and
of right $I^0$-modules.
\ep
\bpf In terms of pseudo compositions, Definition~\ref{D:beta} becomes
\[\beta(X_{(a_0,a_1,\ldots,a_k)})=\begin{cases} 0 &\text{ if }a_0=0\,,\\
X_{(a_0-1,a_1,\ldots,a_k)} &\text{ if }a_0>0\,. \end{cases}\]
Together with~\eqref{E:prodSolB} this shows that $\beta$ is a morphism of right $I^0$-modules.

If $A$ is a commutative algebra, any element $x\in A$ determines a
morphism of bimodules $f:A\to A$ by $f(a)=xa$. Similarly, if $C$ is a cocommutative coalgebra,
any functional $\eta\in C^*$ determines a morphism of bicomodules $f:C\to C$ by
$f=(\eta\ten\id)\circ\Delta$.

Let $C=\Sol{B}$ (a cocommutative coalgebra) and let $\eta:C\to\Q$ be
\[\eta(X_\alpha)=\begin{cases} 0 &\text{ if }\alpha\neq (1)\,,\\
1 & \text{ if }\alpha= (1)\,. \end{cases}\]
It follows from~\eqref{E:coprodSolB}, \eqref{E:prodSolB} and~\eqref{E:defbeta-B} that $\beta=(\eta\ten\id)\circ\Delta$.
Thus, $\beta$ is a morphism of bicomodules.
\epf

We similarly construct a  map  $\pi:\ppp{}\to\ppp{}$ homogeneous of degree $-2$ from the
maps $\pi:\pppn\to\ppp{n-2}$ of Section~\ref{S:intpeaks}.
\bc\label{C:pi-external} The map $\pi:\ppp{}\to\ppp{}$ is a morphism of $\ppp{}$-bicomodules and
of right $\pppint{}$-modules.
\ec
\bpf This follows from  Propositions~\ref{P:beta-external} and~\ref{P:pi}.
\epf

\subsection{Duality. Stembridge's Hopf algebra}\label{S:duality}
In this section we clarify the connection between Stembridge's Hopf algebra of peaks
and the objects discussed in this paper.

The graded dual of $\Sol{A}$ is the Hopf algebra $\QSym$ of quasisymmetric functions.
The dual bases to the bases $\{X_J\}$ and $\{Y_J\}$ of $\Sol{A}$ are, respectively, the {\em monomial}
basis $\{M_J\}$ and Gessel's basis $\{F_J\}$ of $\QSym$.
This is~\cite[Theorem 3.2]{MalReu}, but for the present purposes it may
be taken as the definition of $\QSym$, $M_J$ and $F_J$.

Consider the descents-to-peaks transform $\Theta:\Sol{A}\to\Sol{A}$ of Definition~\ref{E:deftheta}.
The image of this map is $\pppint{}$ (Proposition~\ref{P:theta}).
Let us denote the resulting map $\Sol{A}\onto\pppint{}$ by $j$, and the inclusion
$\pppint{}\inc\Sol{A}$ by $i$, so that $\Theta=i\circ j$.
According to Propositions~\ref{P:hopfpint} and~\ref{P:thetahopf},
the following is a commutative diagram of morphisms of graded Hopf algebras:
\[\smallxymatrix{{\Sol{A}}\ar[rr]^{\Theta}\ar@{->>}[dr]_j & &{\Sol{A}}\\
& {\pppint{}}\ar[ur]_{i} }\]
Dualizing we obtain another diagram of Hopf algebras
 \[\smallxymatrix{{\QSym}\ar[rr]^{\Theta^*}\ar@{->>}[dr]_{i^*} & &{\QSym}\\
& {(\pppint{})^*}\ar[ur]_{j^*} }\]
Now, it follows from~\eqref{E:theta} that the injective map $j^*$ sends
\[(\Pint_F)^*\mapsto \sumsub{J\in\Asets\\F\subseteq J\cup(J+1)}\!\!\!\!2^{1+\#J}\cdot M_J\]
for any $F\in\calFnint$. This says precisely that $j^*$ identifies the Hopf algebra
$(\pppint{})^*$ with the Hopf algebra considered by Stembridge
in~\cite[Theorem 3.1 and Proposition 2.2]{Ste}.

Also, equation~\eqref{E:defpeakinttwo} implies that $i^*$ sends
\[F_J\mapsto(\Pint_{\Lambdaint(J)})^*\,.\]
This says that  $\Theta^*=j^*\circ i^*$ is the map introduced by Stembridge
in~\cite[Theorem 3.1]{Ste}. This provides the explanation announced in Remark~\ref{R:thetadef}.

\br
\be
\item The fact that Stembridge's map $\Theta^*:\QSym\to\QSym$ is a morphism
of Hopf algebras is known from~\cite[Theorem 3.1]{Ste} and~\cite[page 66]{BMSV02}.
In this work we obtained this result
from the fact that its type B analog $\ThetaB$ is (clearly) a morphism of Hopf algebras
(Proposition~\ref{P:thetahopf}).
\item The Hopf algebra structure of $\pppint{}$ (or its dual) has been carefully studied
in~\cite{Sch} and~\cite{Hsi}. One of the important results obtained in these works
is that the elements $\Pint_{\{2,4,\ldots,2k\}}\in\pppint{2k+1}$, $k\geq 0$, freely
generate $\pppint{}$ as an algebra (with respect to the external product).
\item The graded dual of $\ppp{}$ is a comodule algebra over Stembridge's Hopf algebra
$(\pppint{})^*$, and a quotient algebra of $\QSym$. This object has not been considered before.
\ee
\er

\subsection{The action on words}\label{S:words}

In this section we recall the action of $\QS$ on words, discuss its type B analog,
and provide explicit descriptions for the actions of the elements $\Pint_{(n)}$ and
$X^0_{(n)}$. These elements were central to the results of
Sections~\ref{S:principal} and~\ref{S:thetamaps}.

Let $\bfX$ be an arbitrary vector space and $\bfTX=\oplus_{n\geq 0}\bfX^{\ten n}$ its tensor algebra.
The group $\calSn$ acts on $\bfX^{\ten n}$ from the
right by
\[(\bfx_1\ten\cdots\ten \bfx_n)\cdot u=\bfx_{u_1}\ten\cdots\ten \bfx_{u_n}\,.\]
The resulting morphisms of algebras $\Q\calSn\to\End(\bfX^{\ten n})$, $n\geq 0$, may be assembled
into a map
\[\QS\to\End(\bfTX)\]
which by construction sends the internal product of $\QS$ to the composition product of $\End(\bfTX)$.
It is a basic observation of Malvenuto and Reutenauer that, on the other hand, the external
product of $\QS$ corresponds under this map to the {\em convolution} product of
$\End(\bfTX)$~\cite[pages 977-978]{MalReu}. The convolution product is defined from
the canonical Hopf algebra structure of $\bfTX$, for which the elements of $\bfX$ are primitive.

The type B analog of this construction starts from a vector space $\bfY$ endowed
with an involution $\bfy\mapsto\Bar{\bfy}$. There is then a right action of
$\calBn$ on $\bfY^{\ten n}$ given by
\[(\bfy_1\ten\cdots\ten \bfy_n)\cdot w=\bfy_{w_1}\ten\cdots\ten \bfy_{w_n}\,,
\text{ \ where \ }
\bfy_{w_i}=\begin{cases} \bfy_{w_i} & \text{ if }w_i>0\,,\\
\Bar{\bfy}_{-w_i} & \text{ if }w_i<0\,.\end{cases}\]
As before, the resulting map $\QB\to\End(\bfTY)$ sends the internal and external products to the composition
and convolution products of $\End(\bfTY)$, respectively.

Consider the {\em symmetrizer} $\tau:\bfTY\to \bfTY$ defined by
\[\tau(\bfy_1\ten\cdots\ten \bfy_n)= \bfy_1\ten\cdots\ten \bfy_n+
\Bar{\bfy}_n\ten\cdots\ten\Bar{\bfy}_1\,.\]
In other words, $\tau$ acts on $\bfY^{\ten n}$ as the element
$12\ldots n\ +\ \Bar{n}\ldots\Bar{2}\Bar{1}\in\calBn$.
\bp\label{P:action-B} The action of $X^0_{(n)}\in\QB$ on $\bfTY$ is given by
\begin{equation}\label{E:action-B}
(\bfy_1\ten\cdots\ten \bfy_n)\cdot X^0_{(n)}=\tau\Bigl(\cdots\tau\bigl(\tau(\bfy_1)\bfy_2\bigr)
\bfy_3\cdots \bfy_n\Bigr)\,.
\end{equation}
\ep
\bpf We argue by induction on $n$. For $n=1$, $X^0_{(1)}=1+\Bar{1}$, so
\[\bfy\cdot X^0_{(1)}=\bfy+\Bar{\bfy}=\tau(\bfy)\,.\]
Consider the class of signed permutations $X^0_{(n{+}1)}$, that is, the
$(0,n{+}1)$-shuffles~\eqref{E:0nshuffles}. These split into two classes:
those for which the last element is $n{+}1$ and those for which the first element is
$\overline{n{+}1}$. Since in both cases the entries appear in increasing order,
the permutations in the first class are precisely those of the form $\xi\times 1$ and
the permutations in the second class are those of the form
$(\xi\times 1)\cdot\overline{n{+}1}\ldots\Bar{2}\Bar{1}$, with $\xi$ in the class
$X^0_{(n)}$ in both cases. Thus,
\[X^0_{(n{+}1)}=(X^0_{(n)}\times 1)\cdot(12\ldots n{+}1\ +\ \overline{n{+}1}\ldots\Bar{2}\Bar{1})\,.\]
Therefore,
\[(\bfy_1\ten\cdots\ten \bfy_{n{+}1})\cdot X^0_{(n{+}1)}=
\tau\bigl((\bfy_1\ten\cdots\ten \bfy_n)\cdot X^0_{(n)}\ten \bfy_{n{+}1}\bigr)\,.\]
The induction hypothesis then yields the result.
\epf

We go back to the case of an arbitrary vector space $\bfX$.
Consider the {\em Jordan bracket} on $\bfTX$. This is the operator $\bfTX\times \bfTX\to \bfTX$ defined by
\[ \liebrac{\eta,\theta}=\eta\theta+\theta\eta\]
for arbitrary elements $\eta$ and $\theta$ of $\bfTX$.
\bp\label{P:action-P} The action of $\Pint_{(n)}\in\QS$ on $\bfTX$ is given by
\begin{equation}\label{E:action-P}
(\bfx_1\ten\cdots\ten \bfx_n)\cdot \Pint_{(n)}=
\Bigl[\cdots\bigl[ [\bfx_1,\bfx_2],\bfx_3\bigr],\cdots,\bfx_n\Bigr]\,.
\end{equation}
\ep
\bpf  We may view $\bfX$ as endowed with the trivial involution $\Bar{\bfx}=\bfx$ and then
consider the corresponding action of $\QB$ on $\bfTX$. Clearly, this factors through the
action of $\QS$ via the map $\varphi:\QB\to\QS$: $\eta\cdot w=\eta\cdot\varphi(w)$.
Since $\varphi(X^0_{(n)})=2\cdot \Pint_{(n)}$~\eqref{E:shufflesnopeaks}, we have, using~\eqref{E:action-B},
\[(\bfx_1\ten\cdots\ten \bfx_n)\cdot \Pint_{(n)}=
\frac{1}{2}\tau\Bigl(\cdots\tau\bigl(\tau(\bfx_1)\bfx_2\bigr)\bfx_3\cdots \bfx_n\Bigr)\,.\]
Thus, letting $\eta_n=\tau\Bigl(\cdots\tau\bigl(\tau(\bfx_1)\bfx_2\bigr)\bfx_3\cdots \bfx_n\Bigr)$ and
$\theta_n=\Bigl[\cdots\bigl[ [\bfx_1,\bfx_2],\bfx_3\bigr],\cdots,\bfx_n\Bigr]$, we only need to show
 that
 \[\eta_n=2\theta_n\,.\]
For $n=1$ we have
\[\eta_1=\tau(\bfx_1)=\bfx_1+\Bar{\bfx}_1=2\bfx_1=\liebrac{\bfx_1,\bfx_1}=2\theta_1\,,\]
so the result holds. If it does for $n$, then
\[\eta_{n+1}=\tau(\eta_n \bfx_{n+1})=\eta_n \bfx_{n+1}+\bfx_{n+1}\eta_n=
2\theta_n \bfx_{n+1}+ \bfx_{n+1}2\theta_n=2\liebrac{\theta_n, \bfx_{n+1}}=2\theta_{n+1}\,,\]
as needed.
\epf

\begin{rem} Proposition~\ref{P:action-P} is due to Krob, Leclerc and Thibon~\cite[Section 5.3]{KLT}.
A more general operator of the form $\liebrac{\eta,\theta}_q=\eta\theta-q\theta\eta$
is in fact considered in that reference.
\end{rem}



\def\cprime{$'$}
\providecommand{\bysame}{\leavevmode\hbox to3em{\hrulefill}\thinspace}
\providecommand{\MR}{\relax\ifhmode\unskip\space\fi MR }
\providecommand{\MRhref}[2]{%
  \href{http://www.ams.org/mathscinet-getitem?mr=#1}{#2}
}
\providecommand{\href}[2]{#2}

\end{document}